\documentclass{amsart}
\usepackage[latin1]{inputenc}
\usepackage{amssymb}
\usepackage{amsmath}
\usepackage{enumerate}
\usepackage{epsfig,color,caption}

\usepackage[dvipsnames]{xcolor}
\usepackage{colonequals, tikz}
\usepackage{tikz-cd}

\usepackage{color}

\usepackage{longtable,graphics,multirow,ulem}
\usepackage{tikz}
\usepackage{pgf,tikz}
\usetikzlibrary{arrows}
\usepackage[all]{xy}

\newcounter{nootje}
\newtheorem{theorem}{Theorem}[section]

\newtheorem{proposition}[theorem]{Proposition}

\newtheorem{definition}[theorem]{Definition}

\newtheorem{ass}[theorem]{Assumption}

\newtheorem{rem}[theorem]{Remark}

\numberwithin{equation}{section}
\newcommand{\rk}{\mbox{rank}}

\newcommand{\ra}{\rightarrow}

\newcommand{\Fix}{\operatorname{Fix}}
\newcommand{\PGL}{\operatorname{PGL}}

\newcommand{\C }{ \mathbb{C}}

\newcommand{\Z}{\mathbb{Z}}

\usepackage{fancyhdr}
\makeatletter
\def\blfootnote{\xdef\@thefnmark{}\@footnotetext}

\author{Alice Garbagnati}
\address{Dipartimento di Matematica, Univ. Statale di Milano, Milan, Italy}
\email{alice.garbagnati@unimi.it}
\urladdr{ https://sites.google.com/site/alicegarbagnati/}
\author{Cec\'ilia Salgado}
\address{Instituto de Matem\'atica, Univ. Federal do Rio de Janeiro, Rio de Janeiro, Brazil}
\email{salgado@im.ufrj.br}
\urladdr{ http://www.im.ufrj.br/\~{}salgado}

\title[Elliptic fibrations on K3 with non-symplectic involution]{Elliptic fibrations on K3 surfaces with a non-symplectic involution fixing rational curves and a curve of positive genus}
\begin{document}

\keywords{Elliptic fibrations, Rational elliptic surfaces, K3 surfaces, Double covers.\\
The first author is partially supported by FIRB 2012 ``Moduli Spaces and their Applications".\\
The second author is partially supported by Cnpq- Brasil grant nb: 305743/2014-7, FAPERJ E-26/203.205/2016  and the grant Chamada $\#$ 1 of Instituto Serrapilheira.}

\maketitle

\begin{abstract}
In this paper we complete the classification of the elliptic fibrations on K3 surfaces which admit a non-symplectic involution acting trivially on the N\'eron--Severi group. We use the geometric method introduced by Oguiso and moreover we provide a geometric construction of the fibrations classified. If the non-symplectic involution fixes at least one curve of genus 1, we relate all the elliptic fibrations on the K3 surface with either elliptic fibrations or generalized conic bundles on rational elliptic surfaces. This description allows us to write the Weierstrass equations of the elliptic fibrations on the K3 surfaces explicitly and to study their specializations.
\end{abstract}
\section{Introduction}
The purpose of this paper is the classification of elliptic fibrations (with section) on several families of K3 surfaces. These families are characterized by the presence of a non--symplectic involution on their general member.

The families we are interested in were classified by Nikulin, in \cite{Nik non sympl}: let $X$ be a K3 surface on $\C$ and $\iota$ an involution on $X$ which does not preserve the symplectic structure. The fixed locus of $\iota$ can be one of the following:
\begin{itemize}
\item[a)] $\mathrm{Fix}(\iota)= \emptyset$.
\item[b)] $\mathrm{Fix}(\iota)=C \coprod C_i$, where $0\leq i \leq 9$, $C$ is a curve of genus $g\geq 0$ and $C_i$ are rational curves.
\item[c)] $\mathrm{Fix}(\iota)= C \coprod D$, where $C$ and $D$ are both curves of genus 1.
\end{itemize}
The condition that $X$ is a generic among the K3 surfaces admitting a non--symplectic involution with a prescribed fixed locus is equivalent to the condition that $\iota$ acts trivially on its N\'eron-Severi
group.

We classify all elliptic fibrations on $X$ in terms of their
trivial lattice. In particular one obtains that if
$\mathrm{Fix}(\iota)= \emptyset$, then $X$ does not admit elliptic
fibrations with section and if $\mathrm{Fix}(\iota)= C \coprod D$ as in $c)$
then $X$ admits a unique elliptic fibration with section. Hence we concentrate
ourselves on K3 surfaces with a non-symplectic involution whose
fixed locus is as in $b)$. The case $g=0$ was already considered
in \cite{CG}. Here we principally focus on the ones with $g=1$,
i.e., the non-symplectic involution $\iota$ fixes one curve of
genus $1$ and $1\leq k \leq 9$ rational curves, but in Section \ref{sec: higher genus}
we also discuss all the other cases, completing the
classification.

Several papers are devoted to the classifications of elliptic fibrations on K3 surfaces. The most classical are
\cite{Nis}, where a lattice theoretic method is applied, and
\cite{Oguiso} where a more geometric technique is considered. More
recently, the method used in \cite{Nis} is applied, for example, in
\cite{BL} and \cite{wine1} and the one proposed in \cite{Oguiso}
is considered in \cite{Kloosterman}, \cite{CG}, \cite{wine2}. A very deep recent result was obtained in \cite{Ku} where the author
classified the elliptic fibrations on the Kummer surface of a
principally polarized Abelian surface without applying the
previous method, which do not apply well in this situation.

In order to classify the elliptic fibrations on K3 surfaces, we first
produce a list of the possible elliptic fibrations, applying the
techniques developed by Oguiso in \cite{Oguiso}. These techniques
are based on the presence of a non-symplectic involution acting
trivially on the N\'eron--Severi group and here we obtain the
classifications on all the families for which this method applies.
If $g=1$, we construct explicitly the elliptic fibrations listed by means of the geometric constructions presented in our previous paper
\cite{GSal}. These constructions can be considered if the K3
surface is a 2-cover of a rational elliptic surface.

Once one has a classification of the elliptic fibrations on a K3
surface, it is quite natural to ask for the Weierstrass equations
of such fibrations, see e.g. \cite{U}, \cite{KS}, \cite{L},
\cite{BL}. The geometric realization that we provide for the
elliptic fibrations allows us to obtain immediately the
Weierstrass equations by applying an algorithm presented in
\cite{wine2}. Moreover, the knowledge of the equations of the
classified elliptic fibrations allows one to consider
specializations of the elliptic fibrations and thus of the
underlying K3 surfaces. Hence we are able to find some values of
the parameters of the considered families of K3 surfaces for which
the transcendental lattice jumps, and to compute the
new transcendental lattice for these values. In certain cases the
elliptic fibrations specializes because the Mordell--Weil rank
increases and our methods allow us to identify the new sections of the
fibration.

This paper is organized as follows: In Section \ref{sec: general result and admissible fibers}, we state the main Theorem (Theorem \ref{theorem: the classification}) and give a
list of all possible configurations of the trivial lattice of
genus 1 fibrations on the K3 surfaces described above. Then we concentrate, in Sections \ref{sec: RES}, \ref{sec: geometric interpretations}, \ref{sec: equations} and \ref{sec: specializations} on the case $g=1$: Section \ref{sec: RES} is devoted to outlining which rational elliptic surfaces can arise as the quotient $X/ \iota$ mentioned above. Section \ref{sec: geometric interpretations},
contains a realization of the classification of the elliptic
fibrations on $X$  in terms of (generalized) conic bundles
on $X/\iota$. This allows one to compute, in Section \ref{sec: equations}, the Weierstrass
equations of all the elliptic fibrations classified.
In Section \ref{sec: specializations}, we use the equations computed in order to describe
several interesting specializations of the considered K3 surfaces
and of their elliptic fibrations. Section \ref{sec: higher genus} and the Appendix contain the classifications of elliptic fibrations on K3 surfaces which admits a non-symplectic involution acting trivially on the N\'eron--Severi group and fixing one curve of genus greater than 1 and complete the proof of our main Theorem.

\section*{Basic definitions}
In what follows we present the basic key definitions to this text.

\begin{definition} Let $X$ be a smooth projective algebraic surface and $\mathcal{B}$ a smooth projective algebraic curve. An elliptic fibration with base $\mathcal{B}$ on $X$ is a flat morphism $\varepsilon: X \rightarrow B$ such that:
\begin{enumerate}
\item[i)] $\varepsilon^{-1}(t)$ is a smooth curve of genus one for all but finitely many $t \in B$.
\item[ii)] there is a section $\sigma: B \rightarrow X$, i.e., a map such that $\varepsilon \circ \sigma: B \rightarrow B$ is the identity map.
\item[iii)] $\varepsilon^{-1}(t)$ is singular for at least one $t \in B$.
\item[iv)] $\varepsilon$ is relatively minimal, i.e., the fibers of $\varepsilon$ do not admit $(-1)$-curves as components.
\end{enumerate}
\end{definition}

Condition $ii)$ above assures that all but finitely many fibers of
$\varepsilon$ are elliptic curves. Condition $iii)$ rules out
surfaces of product type, i.e., $X \simeq C \times B$, where $C$
is a curve. Finally, note that $iv)$ above means that the surface
is relatively minimal with respect to the fibration, but it does
not imply that $X$ is a minimal surface as $(-1)$-curves are
allowed outside of the fibers of $\varepsilon$. If $\varepsilon$
satisfies $i), ii)$ and $iii)$ but does not satisfy $iv)$ then it
is be called a \textit{non-relatively minimal elliptic fibration}.

We will denote the Mordell--Weil group of $\varepsilon$, i.e. the group of the sections of the fibration, by $MW(\varepsilon)$.
\begin{definition}
A K3 surface is a smooth projective algebraic surface, say $X$, such that:
\begin{itemize}
\item[i)] $q:=h^1(X, \mathcal{O}_X)=0$, i.e., $X$ is regular.
\item[ii)] $K_X \simeq 0$, i.e., the canonical divisor of $X$ is trivial.
\end{itemize}
\end{definition}

If $\varepsilon: X \rightarrow B$ is an elliptic fibration and $X$ is either a K3 surface or a rational surface, then $B\simeq \mathbb{P}^1$.

\begin{definition} Let $X$ be a K3 surface, then $H^{2,0}(X, \mathbb{Z})\simeq \mathbb{Z}\cdot \omega_X$, where $\omega_X$ is a nowhere vanishing symplectic form. An involution $\iota$ on $X$ is called non-symplectic if it does not preserve the symplectic structure on $X$, i.e., $\iota(\omega_X)=-\omega_X$.
\end{definition}

\section{The non-symplectic involution and admissible fibrations }\label{sec: general result and admissible fibers}

The aim of this section is to classify the elliptic fibrations
which appear on surfaces $X$ as in the following assumption:
\begin{ass}\label{cond: X iota}
Let $X$ be a K3 surface and $\iota$ a non--symplectic
involution of $X$ which acts trivially on the N\'eron--Severi
group.
\end{ass}
This assumption is very natural, since this means that the K3
surface is generic in the family of the K3 surfaces with a
non-symplectic involution with a prescribed fixed locus. Nice and
easy examples of these K3 surfaces are provided by double covers
of $\mathbb{P}^2$ branched along a (possible singular or
reducible) sextic or double covers of a rational elliptic surface
such that all the reducible fibers are reduced and contained in
the branch locus.\\

Given a non--symplectic involution $\iota$ we call \emph{special for $\iota$} the curves which are fixed by it. Throughout this note, we  simply call these curves \emph{special} since the dependence on the involution is clear.
\begin{proposition}\label{prop: intersection of special and non-special curves}
Let $(X,\iota)$ be as in Assumption \ref{cond: X iota}.
Let $C$ be
a smooth rational curve in $X$, then either $C$ is special, or it
meets the fixed locus of $\iota$ in exactly two points.

Let $C_1$ and $C_2$ be two smooth rational curves which are not special, then $C_1\cdot C_2\equiv 0\mod 2$
\end{proposition}
\proof This follows immediately from the results from Oguiso and
Kloosterman (see \cite{Oguiso} and \cite{Kloosterman}) and is
due to the facts that the class of each rational curve
is mapped to itself by $\iota$.\endproof
\begin{proposition}\label{lemma: type of fibrations}
Let $(X,\iota)$ be as in Assumption \ref{cond: X iota} and
$\mathcal{E}:X\ra\mathbb{P}^1$ an elliptic fibration on $X$. Then
either:
\begin{itemize}\item[(1)] $\iota$ maps each
fiber of $\mathcal{E}$ to itself or \item[(2)] $\iota$ maps at
least one fiber of $\mathcal{E}$ to another fiber of
$\mathcal{E}$.
\end{itemize}

If $\iota$ is as in (1), then it acts as the identity on the basis
of $\mathcal{E}$ and as the elliptic involution on the fibers.

If $\iota$ is as in (2), then it acts as an involution on the
basis of $\mathcal{E}$ and it preserves two fibers.

The involution $\iota$ is as in (1) if and only if there is a
section of $\mathcal{E}$ which is special curve.
\end{proposition}
\proof Since $\iota$ acts as the identity on the N\'eron--Severi
group, it maps the class of the fiber of $\mathcal{E}$ to itself.
So it maps each fiber either to itself or to another
fiber of the same fibration. In the first case $\iota$ restricts
to the identity of the basis, but since it is not the identity on
$X$ it acts on the fibers. Since $\iota$ is non-symplectic it is
not a translation on each fiber. Hence it is the elliptic
involution on the fibers. The class of the zero section is preserved, so $\iota$ is the elliptic involution and the zero section is a special curve.

If $\iota$ does not preserves each fiber, then it is not the
identity on the basis, and thus it is an involution on the basis,
with two fixed points $p_1$ and $p_2$. Its fixed locus is
necessarily contained in the two fibers over the points $p_1$ and
$p_2$, and so there are no sections among the special curves.
\endproof

\begin{definition}
An elliptic fibration is of type 1 (resp. of type
2) with respect to $\iota$ if it is as in Proposition \ref{lemma: type
of fibrations} case (1) (resp. as in Proposition
\ref{lemma: type of fibrations} case (2)). 
\end{definition}

\begin{proposition}\label{prop: mordell weil}
Let $(X,\iota)$ be as in Assumption \ref{cond: X iota} and
$\mathcal{E}:X\ra\mathbb{P}^1$ an elliptic fibration on $X$. Then
the following hold:
\begin{enumerate}
\item If $\mathcal{E}$ is of type 1 with respect to $\iota$, then
$MW(\mathcal{E})\subset (\Z/2\Z)^2$.
\item If $\mathcal{E}$ is of
type 1 with respect to $\iota$, and there is at least one
non-rational special curve, then $MW(\mathcal{E})\subset \Z/2\Z$.
\item If there is at least one special curve $C$ of genus greater
than 1, then $MW(\mathcal{E})\subset \Z/2\Z$ and in this case if
$MW(\mathcal{E})=\Z/2\Z$, then $C$ is hyperelliptic. \item If
there is at least one special curve of genus greater than 3, then
$MW(\mathcal{E})$ is trivial.
\end{enumerate}
\end{proposition}

\proof The sections of an elliptic fibration are rational curves, hence they have to be mapped to itself by $\iota$. If $\iota$ is the identity on the basis, then each section of the fibration is a fixed curve. In particular each section of an elliptic fibration of type 1 is fixed by $\iota$. If $\mathcal{E}$ is of type 1 with respect to $\iota$, then $\iota$ is the elliptic involution and in particular it fixes the zero section, the (possibly reducible) trisection passing through the 2-torsion points and no other sections. We conclude that if $\mathcal{E}$ is of type 1 with respect to $\iota$, then $MW(\mathcal{E})\subset(\Z/2\Z)^2$. If moreover there is a non rational curve fixed by $\iota$, then it is a component of the trisection passing through the 2-torsion points and thus this trisection is either irreducible (and so $MW(\mathcal{E})=\{0\}$) or it is a bisection (and so $MW(\mathcal{E})=\Z/2\Z$).

We recall that if $\mathcal{E}$ is of type 2 with respect to $\iota$, then the special curves are contained in two fibers, so there are no special curves with genus higher than 1. Hence, if there is at
least one special curve $C$ of genus greater than 1, then $\mathcal{E}$ is of type 1 with respect to $\iota$, and we conclude that either $C$ is the trisection of the 2-torsion points, or the trisection splits into a section and a bisection. In the latter case $C$ is the bisection and by definition it is hyperelliptic.

If $\iota$ fixes one curve of genus higher than 3, then, denoted
by $r$ the rank of the N\'eron--Severi group and by $a$ the length
of the N\'eron--Severi group, it follows by \cite{Nik non sympl}
that $(22-r-a)/2>3$. This implies that $r+a>16$. By \cite{GS sympl
and non sympl}, if $r+a>16$, then $X$ cannot admit a symplectic
involution. On the other hand, if a K3 surface admits an elliptic
fibration with a 2-torsion section, the translation by this
section is a symplectic involution on the K3 surface. We conclude
that if $r+a>16$, there are no elliptic fibrations on $X$ with a
2-torsion section.
\endproof

\begin{theorem}\label{theorem: the classification}\label{prop: when iota is trivial on the basis}
Let $X$ be a K3 surface which admits at least one elliptic
fibration and $(X,\iota)$ as in Assumption \ref{cond: X iota}. The
following hold:
\begin{itemize}
\item[i)] $\Fix_{\iota}(X)\neq \emptyset$;
\item[ii)] If $\Fix_{\iota}(X)=C\coprod D$ with $g(C)=g(D)=1$, then there is a unique elliptic fibration on $X$, which is $\varphi_{|C|}:X\ra\mathbb{P}^1$ and it is of type 2;
\item[iii)] If $\mathrm{Fix}(\iota)=C \coprod C_i$, $g(C_i)=0$, and $g(C)>1$ then all the fibrations on $X$ are of type 1 and they are given in Proposition \ref{prop: higher genus}
\item[iv)] If $\mathrm{Fix}(\iota)=C \coprod C_i$, $g(C_i)=0$ and $g(C)=1$ then there exists one fibration of type 2, $\varphi_{|C|}:X\ra\mathbb{P}^1$, and all the other fibrations are of type 1. All the fibrations on $X$ are given in Proposition \ref{prop: possible fibrations}.
\item[v)] If $\mathrm{Fix}(\iota)=C \coprod C_i$, $g(C_i)=g(C)=0$ and then there exists both fibrations of type 2 and of type 1. All the fibrations on $X$ are given in \cite{CG}.
\end{itemize}
\end{theorem}
\proof If $\mathcal{E}$ is an elliptic fibration of type 1 with
respect to $\iota$, then the zero section is a special curve
$\mathcal{E}$. If $\mathcal{E}$ is an elliptic fibration of type
2, then it preserves two fibers of the fibration and the zero
section. So $\iota$ cannot be a translation on these fibers and
thus $\Fix_{\iota}(X)\neq \emptyset$.

Let us assume that $\iota$ fixes at least one genus 1 curve $E$.
Then $\varphi_{|E|}:X\ra\mathbb{P}^1$ is a genus 1 fibration and
$\iota$ fixes at least one of the fibers. If $\iota$ would act as
the identity on the base of the fibration it cannot act as the
identity also on the fibers of the fibration (otherwise it is the
identity). So, if $\iota$ is the identity on the fiber $E$ of
$\varphi_{|E|}:X\ra\mathbb{P}^1$, it is not the identity on the
base, i.e. $\varphi_{|E|}$ is of type 2 with respect to $\iota$.
In particular it is an involution of the basis, which fixes two
points on the basis, and thus it preserves two fibers of
$\varphi_{|E|}:X\ra \mathbb{P}^1$. The special curve $E$ is one of
these fibers.

Assume that $X$ admits other genus 1 fibrations, which are not
$\varphi_{|E|}$. Denote by $\mathcal{E}:X\ra \mathbb{P}^1$ one of
these. The special curve $E$ cannot be a fiber of $\mathcal{E}$ as
otherwise $\mathcal{E}$ would coincide with $\varphi_{|E|}$. So $E$ is a horizontal curve and $\mathcal{E}$ is of type 1 with resect to $\iota$. Since $E$ is not a rational curve, it is neither a section nor an irreducible component of a reducible fiber. Thus it is a multisection passing through to 2-torsion points of the fibers of $\mathcal{E}$, i.e. it is either a trisection or a
bisection. In both the cases there can not be another genus 1
curve in the fixed locus.

If there is a special curve of genus bigger than 1, then it is not contained in a fiber and thus there are no fibrations of type 2 with respect to $\iota$.

Thanks to the results in \cite{Oguiso}, \cite{Kloosterman} and
\cite{CG}, to conclude the proof it remains to classify the
elliptic fibrations which appear in the cases $\mathrm{Fix}(\iota)=C
\coprod C_i$ with $g(C)\geq 1$. This is done in Propositions
\ref{prop: possible fibrations} and \ref{prop: higher genus}.\endproof
We observe that most of elliptic fibrations that we
are looking for are of type 1. For these fibrations the Mordell--Weil group is extremely simple, so that one has to classify principally the reducible fibers which can appear. This is the purpose of the following proposition, where we reformulate the results by Oguiso, see \cite{Oguiso}, to deal with surfaces in our setting.

\begin{proposition}\label{prop: possible fibers}
Let $(X,\iota)$ be as in Assumption \ref{cond: X iota} and
$\mathcal{E}:X\ra\mathbb{P}^1$ an elliptic fibration on $X$ of
type 1 with respect to $\iota$. Then the reducible fibers which
appear in $\mathcal{E}$ are the ones contained in the table, where
the number $s$ of special rational curves that they contain and
the number $c$ of non-trivial components that they contain  are
given.
$$
\begin{array}{|c|c|c|c|c|c|c}
\hline
fiber&s&c& Dynkin\\
\hline
I_2&0&1&A_1\\
\hline
I_{2n}^*&n+1&2n+4&D_{2n+4}\\
\hline
III^*&3&7&E_7\\
\hline
II^*&4&8&E_8\\
\hline
\end{array}
$$ \end{proposition}
\proof This is a trivial consequence of Propositions \ref{prop:
intersection of special and non-special curves} and \ref{prop:
mordell weil}. In particular, by the assumptions,
$\iota$ fixes the zero section of the fibration.

To illustrate the kind of arguments, we show why fibers of type
$IV^*$ are not allowed. Suppose that there exists such a fiber. We
denote by $\Theta_i$ its components: $\Theta_0$ is the components
meeting the zero section. If $i,j\in\{0,1,\ldots, 4\}$, then
$\Theta_i\cdot \Theta_j=1$ if and only if $|i-j|=1$ and
$\Theta_i\cdot \Theta_j=0$ otherwise. The component $\Theta_5$
meets only the component $\Theta_2$ in one point and the component
$\Theta_6$ in another point. The component $\Theta_6$ meets only
the component $\Theta_5$. Since $\iota$ acts trivially one the
N\'eron--Severi group, each component $\Theta_i$
is sent to itself by $\iota$. In particular
the intersection points between $\Theta_2$ and $\Theta_i$, for
$i=1,3,5$ are fixed points. But then $\Theta_2$ is a fixed curve
(the involution $\iota$ acts on the rational curve $\Theta_2$ with
3 fixed points). The fixed locus of $\iota$ is smooth, so
$\Theta_1$, $\Theta_3$ and $\Theta_5$ are not fixed. The
intersection point between $\Theta_1$ and $\Theta_0$ is a fixed
point and it is a singular points of the fiber $IV^*$. Neither a section nor the trisection of the 2-torsion points pass through this singular point of $IV^*$. A non-symplectic
involution on a K3 surface cannot admit an isolated fixed point,
so there is a curve passing through the intersection point between
$\Theta_1$ and $\Theta_0$. So $\Theta_0$ is a
fixed curve. Analogously $\Theta_4$ and $\Theta_6$ are fixed. But
these curves are the unique simple components of the fiber of type
$IV^*$ and we know that $\iota$ fixes at least on section. This means that there exist two special curves which intersect, namely a section of the fibration and a component of a fiber, which is impossible by the smoothness of the fixed locus of $\iota$.

\begin{proposition}\label{prop: K3 with non symple involution g=1}{\rm (See \cite{Nik non sympl})}
Let $(X,\iota)$ be as in Assumption \ref{cond: X iota}. Let us
assume that $\iota$ fixes a curve of genus 1 and other $k$
rational curves. Then:\begin{itemize}\item[i)] the N\'eron--Severi
group of $X$ has rank $r=10+k$ and its discriminant is
$(\Z/2\Z)^a$, where $a=20-r(=10-k)$;\item[ii)] the pair $(r,a)$
determines $NS(X)$ if $r\neq 14$ and $r\neq 18$. If $r=14$ or
$r=18$ there are two different possibilities for $NS(X)$, which
depend on the values of $\delta\in\{0,1\}$. If the discriminant
form of $NS(X)$ takes values in $\Z$, then $\delta=0$, otherwise
$\delta=1$; \item [iii)]The triple $(r,a,\delta)$ determines
uniquely $NS(X)$.
\end{itemize}
\end{proposition}

\begin{proposition}\label{prop: possible fibrations}
Let $(X,\iota)$ be as in Proposition \ref{prop: K3 with non symple
involution  g=1} and $\mathcal{E}:X\ra\mathbb{P}^1$ an elliptic
fibration on $X$. Then there are two possibilities:
\begin{itemize}\item[i)] $\mathcal{E}$ is the unique fibration of type 2 and the configuration for the reducible fibers
appears as first lines in Table \ref{table of admissible fibers
g=1}, \item[ii)] $\mathcal{E}$ is of type 1. The admissible
configurations of the reducible fibers are listed in Table
\ref{table of admissible fibers g=1}, in all the other
lines.\end{itemize}

\begin{longtable}{l}\caption{}\label{table of admissible fibers g=1}\\ \hline k=9, r=19, a=1\\
$
\begin{array}{|c|c|c|c|c|c|}
\hline
n^o&\mbox{trivial lattice}&17=\sum c_i+\rk(MW)&9=k=\sum s_i+\#\mbox{sections}&MW(\mathcal{E})\\
\hline
9.1&U\oplus A_{17}&17+0&9+0&\Z/3\Z\\
\hline
9.2&U\oplus E_8\oplus E_8\oplus A_1&8+8+1+0& 4+4+0+1&\{1\}\\
\hline
9.3&U\oplus E_7\oplus D_{10}&7+10+0& 3+5+2&\Z/2\Z\\
\hline
9.4&U\oplus D_{16}\oplus A_1&16+1+0& 7+0+2&\Z/2\Z\\
\hline
\end{array}$\\
 $ $\\
\hline
\mbox{k=8, r=18,\ a=2, }$\delta=0$\\
$\begin{array}{|c|c|c|c|c|c|} \hline
n^o&\mbox{trivial lattice}&16=\sum c_i+\rk(MW)&8=k=\sum s_i+\#\mbox{sections}&MW(\mathcal{E})\\
\hline
8.1&U\oplus A_{15}&15+1&8+0&\Z\times \Z/2\Z\\
\hline
8.2&U\oplus E_8\oplus D_{8}&8+8+0& 4+3+1&\{1\}\\
\hline
8.3&U\oplus E_7^{ \oplus 2}\oplus A_1^{\oplus 2}&7+7+1+1+0&3+3+0+0+2&\Z/2\Z\\
\hline
8.4&U\oplus D_{16}&16+0&7+1&\{1\}\\
\hline
8.5&U\oplus D_{12}\oplus D_4&12+4+0&5+1+2&\Z/2\Z\\
\hline
8.6&U\oplus D_{8}\oplus D_8&8+8+0&3+3+2&\Z/2\Z\\
\hline
\end{array}$\\
$ $\\

\hline
\mbox{k=8, r=18,\ a=2, }$\delta=1$\\
$\begin{array}{|c|c|c|c|c|c|} \hline
n^o&\mbox{trivial lattice}&16=\sum c_i+\rk(MW)&8=k=\sum s_i+\#\mbox{sections}&MW(\mathcal{E})\\
\hline
8.1&U\oplus A_{15}&15+1&8+0&\Z\\
\hline
8.2&U\oplus E_8\oplus E_7\oplus A_1&8+7+1+0& 4+3+0+1&\{1\}\\
\hline
8.3&U\oplus E_7\oplus D_8\oplus A_1&7+8+1+0&3+3+0+2&\Z/2\Z\\
\hline
8.4&U\oplus D_{14}\oplus A_1^{\oplus 2}&14+1+1+0&6+0+0+2&\Z/2\Z\\
\hline
8.5&U\oplus D_{10}\oplus D_6&10+6+0&4+2+2&\Z/2\Z\\
\hline
\end{array}$\\
$ $\\

\hline

\mbox{k=7, r=17,\ a=3}\\
$\begin{array}{|c|c|c|c|c|c|} \hline
n^o&\mbox{trivial lattice}&15=\sum c_i+\rk(MW)&7=k=\sum s_i+\#\mbox{sections}&MW(\mathcal{E})\\
\hline
7.1&U\oplus A_{13}&13+2&7+0&\left(\Z\right)^2\\
\hline
7.2&U\oplus E_8\oplus D_6\oplus A_1&8+6+1+0& 4+2+0+1&\{1\}\\
\hline
7.3&U\oplus E_7\oplus D_6\oplus A_1^{\oplus 2}&7+6+1+1+0& 3+2+0+0+2&\Z/2\Z\\
\hline
7.4&U\oplus E_7\oplus D_8&7+8+0& 3+3+1&\{1\}\\
\hline
7.5&U\oplus E_7^{ \oplus 2}\oplus A_1&7+7+1+0& 3+3+0+1&\{1\}\\
\hline
7.6&U\oplus D_{14}\oplus A_1&14+1+0& 6+0+1&\{1\}\\
\hline
7.7&U\oplus D_{12}\oplus A_1^{\oplus 3}&12+1+1+1+0& 5+0+0+0+2&\Z/2\Z\\
\hline
7.8&U\oplus D_{10}\oplus D_4\oplus A_1&10+4+1+0& 4+1+0+2&\Z/2\Z\\
\hline
7.9&U\oplus D_{8}\oplus D_6\oplus A_1&8+6+1+0& 3+2+0+2&\Z/2\Z\\
\hline
\end{array}$\\
 $ $\\
\hline
\mbox{k=6, r=16,\ a=4}\\
$\begin{array}{|c|c|c|c|c|c|} \hline
n^o&\mbox{trivial lattice}&14=\sum c_i+\rk(MW)&6=k=\sum s_i+\#\mbox{sections}&MW(\mathcal{E})\\
\hline
6.1&U\oplus A_{11}&11+3&6+0&\left(\Z\right)^3\\
\hline
6.2&U\oplus E_8\oplus D_4\oplus A_1^{\oplus 2}&8+4+1+1+0& 4+1+0+0+1&\{1\}\\
\hline
6.3&U\oplus E_7\oplus D_6\oplus A_1&7+6+1+0& 3+2+0+1&\{1\}\\
\hline
6.4&U\oplus E_7\oplus D_4\oplus A_1^{\oplus 3}&7+4+1+1+1+0& 3+1+0+0+0+2&\Z/2\Z\\
\hline
6.5&U\oplus D_{12}\oplus A_1^{\oplus 2}&12+1+1+0&5+0+0+1&\{1\}\\
\hline
6.6&U\oplus D_{10}\oplus D_4&10+4+0& 4+1+1&\{1\}\\
\hline
6.7&U\oplus D_{10}\oplus A_1^{\oplus 4}&10+1+1+1+1+0& 4+0+0+0+0+2&\Z/2\Z\\
\hline
6.8&U\oplus D_{8}\oplus D_6&8+6+0& 3+2+0+0+1&\{1\}\\
\hline
6.9&U\oplus D_{8}\oplus D_4\oplus A_1^{\oplus 2}&8+4+1+1+0& 3+1+0+0+2&\Z/2\Z\\
\hline
6.10&U\oplus D_{6}^{\oplus 2}\oplus A_1^{ \oplus 2}&6+6+1+1+0& 2+2+0+0+2&\Z/2\Z\\
\hline
\end{array}$\\
$ $\\
\hline
\mbox{k=5, r=15,\ a=5}\\
$\begin{array}{|c|c|c|c|c|c|} \hline
n^o&\mbox{trivial lattice}&13=\sum c_i+\rk(MW)&5=k=\sum s_i+\#\mbox{sections}&MW(\mathcal{E})\\
\hline
5.1&U\oplus A_{9}&9+4&5+0&\left(\Z\right)^4\\
\hline
5.2&U\oplus E_8\oplus A_1^{\oplus 5}&8+1+1+1+1+1+0& 4+0+0+0+0+1&\{1\}\\
\hline
5.3&U\oplus E_7\oplus D_4\oplus A_1^{\oplus 2}&7+4+1+1+0& 3+1+0+0+1&\{1\}\\
\hline
5.4&U\oplus E_7\oplus A_1^{\oplus 6}&7+1+1+1+1+1+1+0& 3+0+0+0+0+0+0+2&\Z/2\Z\\
\hline
5.5&U\oplus D_{10}\oplus A_1^{\oplus 3}&10+1+1+1+0& 4+0+0+0+1&\{1\}\\
\hline
5.6&U\oplus D_{8}\oplus D_4\oplus A_1&8+4+1+0& 3+1+0+1&\{1\}\\
\hline
5.7&U\oplus D_{8}\oplus A_1^{\oplus 5}&8+1+1+1+1+1+0& 3+0+0+0+0+0+2&\Z/2\Z\\
\hline
5.8&U\oplus D_{6}^{ \oplus 2}\oplus A_1&6+6+1+0& 2+2+0+1&\{1\}\\
\hline
5.9&U\oplus D_{6}\oplus D_4\oplus A_1^{\oplus 3}&6+4+1+1+1+0& 2+1+0+0+0+1&\Z/2\Z\\
\hline
\end{array}$\\
$ $\\
\hline
\mbox{k=4, r=14, a=6 } $\delta=0$\\
$\begin{array}{|c|c|c|c|c|c|} \hline
n^o&\mbox{trivial lattice}&12=\sum c_i+\rk(MW)&4=k=\sum s_i+\#\mbox{sections}&MW(\mathcal{E})\\
\hline
4.1&U\oplus E_{6}&6+6&4+0&\left(\Z\right)^6\\
\hline
4.2&U\oplus D_6\oplus A_1^{\oplus 6}&6+1+1+1+1+1+1+0& 2+0+0+0+0+0+0+2&\Z/2\Z\\
\hline
4.3&U\oplus D_4^{\oplus 3}&4+4+4+0& 1+1+1+1&\{1\}\\
\hline
\end{array}$\\
 $ $\\
\hline
\mbox{k=4, r=14, a=6 }$\delta=1$ \\
$\begin{array}{|c|c|c|c|c|c|} \hline
n^o&\mbox{trivial lattice}&12=\sum c_i+\rk(MW)&4=k=\sum s_i+\#\mbox{sections}&MW(\mathcal{E})\\
\hline
4.1&U\oplus A_{7}&7+5&4+0&\left(\Z\right)^5\\
\hline
4.2&U\oplus E_7\oplus A_1^{\oplus 5}&7+1+1+1+1+1+0& 3+0+0+0+0+0+1&\{1\}\\
\hline
4.3&U\oplus D_8\oplus A_1^{\oplus 4}&8+1+1+1+1+0& 3+0+0+0+0+1&\{1\}\\
\hline
4.4&U\oplus D_6\oplus D_4\oplus A_1^{\oplus 2}&6+4+1+1+0& 2+1+0+0+1&\{1\}\\
\hline
4.5&U\oplus D_6\oplus A_1^{\oplus 6}&6+1+1+1+1+1+1+0& 2+0+0+0+0+0+0+2&\Z/2\Z\\
\hline
4.6&U\oplus D_4^{\oplus 2}\oplus A_1^{\oplus 4}&4+4+1+1+1+1+0& 1+1+0+0+0+0+2&\Z/2\Z\\
\hline
\end{array}$\\
 $ $\\
\hline
\mbox{k=3, r=13, a=7}\\
$\begin{array}{|c|c|c|c|c|c|} \hline
n^o&\mbox{trivial lattice}&11=\sum c_i+\rk(MW)&3=k=\sum s_i+\#\mbox{sections}&MW(\mathcal{E})\\
\hline
3.1&U\oplus A_{5}&5+6&3+0&\left(\Z\right)^6\\
\hline
3.2&U\oplus D_6\oplus A_1^{\oplus 5}&6+1+1+1+1+1+0& 2+0+0+0+0+0+1&\{1\}\\
\hline
3.3&U\oplus D_4^{\oplus 2}\oplus A_1^{\oplus 3}&4+4+1+1+1+0& 1+1+0+0+0+1&\{1\}\\
\hline
3.4&U\oplus D_4\oplus A_1^{\oplus 7}&4+1+1+1+1+1+1+1+0& 1+0+0+0+0+0+0+0+2&\Z/2\Z\\
\hline
\end{array}$\\
$ $\\
\hline
\mbox{k=2, r=12, a=8}\\
$\begin{array}{|c|c|c|c|c|c|} \hline
n^o&\mbox{trivial lattice}&10=\sum c_i+\rk(MW)&2=k=\sum s_i+\#\mbox{sections}&MW(\mathcal{E})\\
\hline
2.1&U\oplus A_{3}&3+7&2+0&\left(\Z\right)^7\\
\hline
2.2&U\oplus D_4\oplus A_1^{\oplus 6}&4+1+1+1+1+1+1+0& 1+0+0+0+0+0+0+1&\{1\}\\
\hline
2.3&U\oplus A_1^{\oplus 10}&10+0& 0+\ldots+0 +2&\Z/2\Z\\
\hline
\end{array}$\\
 $ $\\
\hline
\mbox{k=1, r=11, a=9}\\
$\begin{array}{|c|c|c|c|c|c|} \hline
n^o&\mbox{trivial lattice}&9=\sum c_i+\rk(MW)&1=k=\sum s_i+\#\mbox{sections}&MW(\mathcal{E})\\
\hline
1.1&U\oplus A_{1}&1+8&1+0&\left(\Z\right)^8\\
\hline
1.2&U\oplus A_1^9&9+0& 0+1&\{1\}\\
\hline
\end{array}$
\end{longtable}
\end{proposition}

\proof Assume first that none of the special curves are
sections for $\mathcal{E}$. Then $\mathcal{E}$ is of type 2 and
the quotient $X/\iota$ is a rational elliptic surface $R$ equipped
with the elliptic fibration $\mathcal{E}_R$. The elliptic
fibration $\mathcal{E}$ is induced by $\mathcal{E}_R$ by a base
change of order 2, as proved in \cite[Theorem 4.2]{GSal}, see also
\cite{Z}. Moreover, if $\iota$ acts trivially on the
N\'eron--Severi group, then the elliptic fibration on $R$ has
exactly one reducible fiber which is contained in the branch locus
of the cover $X\ra R$. The reducible fiber is determined by $k$
and $\delta$, and it is always $I_k$, except for $k=4$, when it
can be also $IV$. This last statement is proved in the beginning
of Section \ref{sec: RES}.

If there is a special curve which is a section of $\mathcal{E}$,
then $\mathcal{E}$ is of type 1 and in order to produce the list
one has to makes the following observations:
\begin{itemize}
\item By Proposition \ref{prop: mordell weil}, $MW(\mathcal{E})\subset \Z/2\Z$.
\item Since $MW(\mathcal{E})\subset \Z/2\Z$, the trivial lattice of the fibration has rank $r=\rk(NS(X))$. So the sum of the non-trivial components of the reducible fibers has to be $r-2$.
\item If $MW(\mathcal{E})=0$, then the discriminant of the trivial lattice coincides with the discriminant of $NS(X)$, so it is $(\Z/2\Z)^a$. In this case there is a unique section (the zero section), and it is a special curve. So the sum of the special curves contained in the reducible fibers has to be $k-1$.
\item If $MW(\mathcal{E})=\Z/2\Z$, then the trivial lattice is a sublattice with index 2 in $NS(X)$, so its discriminant is $(\Z/2\Z)^{a+2}$. In this case there are two sections, and they are both special curves. So the sum of the special curves contained in the reducible fibers has to be $k-2$.
\end{itemize}
The list given in the proposition is the list of all the trivial
lattices which satisfy these conditions. In order to conclude that
every elliptic fibrations in the list really occurs we explicitly
construct all of them in the next sections. Alternatively, one can
observe that the fibrations of type 1 have Mordell Weil group with
rank 0, so they exist if and only if they appear in the list
given in \cite{Shim}. The existence of the ones of type 2 follows
by the existence of the associated rational elliptic surface as in Section \ref{sec: RES}.
\endproof

\begin{rem}{\rm A different way to obtain Table \ref{table of admissible fibers  g=1}, is to apply the so called Nishiyama method, \cite{Nis}. In order to apply the method one has to consider a negative definite lattice $T$ with the same discriminant group and form of the transcendental lattice of $X$, $T_X$, and such that $\rk(T)=\rk(T_X)+4$. For the surfaces considered in Proposition \ref{prop: possible fibrations}, if $\delta=1$, $T\simeq E_7\oplus A_1^{9-k}$; if $\delta=0$ and $k=8$, $T\simeq D_8$; if $\delta=0$ and $k=4$, $T\simeq D_4\oplus D_4\oplus D_4$}.\end{rem}

\section{Classification of the admissible rational elliptic surfaces}\label{sec: RES}
We study the rational elliptic surface $R$ which
appears when taking the quotient  $X/\iota$, giving a realizations of the first line of the tables in Proposition \ref{prop: possible
fibrations}.

Given $(X,\iota)$ as in Proposition \ref{prop: possible fibrations} let $\mathcal{E}$ be the fibration of type 2 on $X$. We denote by  $R$ the surface obtained after blowing down $(-1)$-curves that are
components of the fibers on the non relatively minimal fibration induced on $X/ \iota$ by $\mathcal{E}$.

Vice versa, given a rational elliptic surface $\mathcal{E}_R:R\ra \mathbb{P}^1$, let $f:\mathbb{P}^1\ra\mathbb{P}^1$ be a  2:1 map branched over $0$ and $\infty$.
We consider the associated base change. This produces another elliptic fibration $\mathcal{E}_X:X\ra\mathbb{P}^1$, where $X$ is the eventual desingularization of the fiber product $R \times_{\mathbb{P}^1}\mathbb{P}^1$. The surface $X$ naturally comes with an involution $\iota$, such that $X/\iota$ is birational to $R$. We now require the following
\begin{itemize}
\item[a)] $X$ is a K3 surface; \item [b)] $\iota$ fixes one curve
of genus 1 and $1\leq k\leq 9$ rational curves; \item[c)] $\iota$
acts trivially on $NS(X)$.\end{itemize} Condition $a)$ implies
that the fibers over $0$ and $\infty$ of $\mathcal{E}_R$ are
reduced, so they can be $I_n$, $n\geq 0$ or $II$ or $III$ or $IV$
(see \cite{SS}).

Condition $b)$ implies that exactly one among the fibers over $0$ and $\infty$ of $\mathcal{E}_R$ is of type $I_0$ (see \cite{GSal}). Just to fix the notation we assume that the fiber over $0$ is of type $I_0$.

Condition $c)$ implies that $\mathcal{E}_R$ has no reducible fibers except possibly the fiber over $0$ and the fiber over $\infty$ (see \cite{GSal}). Moreover condition $c)$ implies that $\mathcal{E}_R$ has no fibers of type $II$ and $III$ over $0$ and $\infty$. Indeed let us assume that the fiber over $\infty$ is a fiber of type $II$, i.e. a cuspidal curve. Then the fiber over $\infty$ of $\mathcal{E}_X$ is a fiber of type $IV$ which has 3 components: one of them is the double cover of the strict transform of the fiber of type $II$ and it is fixed by $\iota$. The others are two curves switched by $\iota$ (which in fact are mapped to the same curve on the blow up of $R$ in the singular point of the cuspidal fiber). But if $\iota$ switches the components of a reducible fiber, it can not be the identity on $NS(X)$. Similarly, if the fiber of $\mathcal{E}_R$ over $\infty$ is of type $III$, the corresponding fiber over $\mathcal{E}_X$ is of type $I_0^*$. Two of the simple components of the fiber $I_0^*$ are fixed by $\iota$, the multiple component is preserved, but not fixed by $\iota$, an so the other two simple components are switched by $\iota$. Therefore we can exclude also this case.

The admissible rational elliptic surfaces are described in Table
\ref{list: admissible RES} where we give information both on $R$
and on $(X,\iota)$, namely the reducible fibers of $R$; the
Mordell Weil group of $\mathcal{E}_R$; $(r,a,\delta)$, which
determines $NS(X)$; $k$, the number of rational curves fixed by
$\iota$.
\begin{align}\label{list: admissible RES}
\mbox{ }\end{align}
\vspace{-20pt}
\begin{eqnarray*}
\begin{array}{|c||c|}
\hline
\mbox{ Surface }R&\mbox{ Surface }X\\
\hline
\begin{array}{c|c|c}
\mbox{ Branch fibers}&\mbox{ Others singular fibers }&MW\\
\hline
I_0+I_n,\ 1\leq n\leq 8&\mbox{ irreducible }&\Z^{9-n} \\
\hline
I_0+I_9&\mbox{ irreducible }&\Z/3\Z \\
\hline
I_0+I_8&\mbox{ irreducible }&\Z\times \Z/2\Z \\
\hline
I_0+IV&\mbox{ irreducible }&\Z^6 \\

\end{array}
&
\begin{array}{c|c}(r,a,\delta)&k\\
\hline
(10+n,10-n,1)&n\\
\hline
(19,1,1)&9\\
\hline
(18,2,0)&8\\
\hline
(14,6,0)&4\\
\end{array}\\
\hline
\end{array}
\end{eqnarray*}

Notice that, by the list presented in \cite{P}, the rational elliptic surfaces with a unique reducible fiber which is of type $I_n$ have no torsion sections with only two exceptions, $n=8,9$. If $n=8$ there exist two different families of rational elliptic surfaces, the Mordell--Weil group of one of them is torsion free while the other has a 2-torsion section. If $n=9$, then there is a unique rational surface with a reducible fiber $I_9$, it is extremal and the Mordell--Weil group is necessarily $\Z/3\Z$.

\subsection{Equations of the pencil of cubics}
The rational elliptic surface associated to the involution which
fixes 9 rational curves and one curve of genus 1 is the unique
extremal rational surface with a fiber of type $I_9$ and three
fibers of type $I_1$. It is associated to a known pencil of cubics
which is generated, for example, by the reducible cubics $xyz$ and
by the smooth cubic $\mathcal{C}_9:=\{x^2y+y^2z+z^2x=0\}$. In the
following we will call $l_1$ the line $x=0$, $l_2$ the line $y=0$
and $l_3$ the line $z=0$. Moreover we will denote by $P_{ij}$ the
point $l_i\cap l_j$. The cubic $\mathcal{C}_9$ passes through the
points $P_{ij}$ and in $P_{12}$ it is tangent to $l_1$, in
$P_{13}$ it is tangent to $l_3$ and in $P_{23}$ it is tangent to
$l_2$. Hence in each point $P_{ij}$ the cubic $\mathcal{C}_9$
intersects the cubic $l_1\cup l_2\cup l_3$ with multiplicity 3 and
all the base points of the pencil $x^2y+y^2z+z^2x+\mu xyz$ are
$P_{ij}$ and points infinitely near to $P_{ij}$. After blowing up
the base points one obtains a reducible fiber of type $I_9$ over
$\mu=\infty$ and no other reducible fiber.

The rational elliptic surface with one fiber
of type $I_8$ and no torsion is obtained by deforming
the previous example. Indeed, we can obtain a fiber of type $I_8$
over $\mu=\infty$ if we separate one of the base points
infinitely near to $P_{ij}$. This is equivalent to require
that in one point $P_{ij}$ the cubic $\mathcal{C}_9$ deforms to
cubics that still pass through $P_{ij}$, but are not tangent
to a line between $l_i$. We can assume that in $P_{13}=(0:1:0)$
the deformation of $\mathcal{C}_9$ is not tangent to $l_3$. This gives the pencil
$\mathcal{C}_8:=\{x^2y+y^2z+z^2x+a_9xy^2=0\}$.

Proceeding by iterations of the above, we obtain that the pencil
\begin{equation}\label{eq: pencil of
Ih}\mathcal{P}_1:=\{x^2y+y^2z+z^2x+a_9xy^2+a_8x^2z+a_7yz^2+a_6z^3+a_5y^3+a_4x^3+\mu(xyz+a_3z^3+a_2y^3)\}.\end{equation}
corresponds to a rational elliptic surface with one fiber of type
$I_1$ in $\mu=\infty$ and no reducible fibers.

\begin{proposition}\label{prop: equations of RES}
Let $\mathcal{P}_k$ be the pencil of cubics obtained
by choosing $a_2=\cdots =a_k=0$ for $2\leq k\leq 9$ in \eqref{eq: pencil of Ih}. For a generic choice of the $a_i$'s, $\mathcal{P}_k$ corresponds to an elliptic fibration with a
fiber of type $I_k$ over $\mu=\infty$ and no other reducible
fibers.

For a generic choice of $b$, the pencil of cubics $x^2z+y^2z+y^2x+bxz^2+\mu(xyz)$ corresponds to an elliptic fibration with a fiber of type $I_8$ over $\mu=\infty$, no other reducible
fibers and a 2-torsion section.

For a generic choice of the $c_i$'s, the pencil of cubics $z^3+c_1xy^2+c_2x^2z+c_3xyz+c_4y^2z+c_5xz^2+(-1-c_1-c_2-c_3-c_4-c_5)yz^2+\mu xy(x-y)$ corresponds to an elliptic fibration with a fiber of type $IV$ over $\mu=\infty$ and no other reducible
fibers.
\end{proposition}
\proof This is straightforward once one considers the base points
of each pencil $\mathcal{P}_k$.

Alternatively, one can compute the
Weierstrass equation of the elliptic fibration induced by
each of the pencils in the statement and consider the discriminant.\endproof

\section{Geometric construction of type 1 elliptic fibrations}\label{sec: geometric interpretations} The aim of this
section is to provide geometric realizations of the elliptic
fibrations of type 1 on $X$ which are listed in Proposition \ref{prop: possible fibrations}.

This is done by considering linear systems on the quotient surface $X/\iota$, which is a rational surface, denoted from now on by $\widetilde{R}$.

\begin{definition}\label{defi: generalized conic bundle}{\rm (See \cite[Definition 3.3]{GSal})}
A generalized conic bundle on $\widetilde{R}$ is a \textit{nef}
class $D$ in $\mathrm{NS}(\widetilde{R})$ such that
\begin{itemize}
\item[i)]  $D\cdot(-K_{\widetilde{R}})= 2$.
\item[ii)] $D^2=0$.
\end{itemize}
\end{definition}

We observe that $\widetilde{R}$ is a blow up of the rational
elliptic surface $R$ and the previous definition generalizes the
standard definition of conic bundles on $R$. Note that $\tilde{R}$
is endowed with a non-relatively minimal elliptic fibration
induced by the elliptic fibration on $R$. Since $R$ is a rational
elliptic surface it comes with a map $R\ra\mathbb{P}^2$ given by
the blow up of the base points of the pencil of cubics described
in Proposition \ref{prop: equations of RES}. The map
$\widetilde{R}\ra R$ contracts $(-1)$-curves contained in fibers.
Hence we have a contraction map $\widetilde{R}\ra\mathbb{P}^2$.
Some of the generalized conic bundles remain base point free
systems on $R$, and define standard conic bundles on $R$. All the
generalized conic bundles are mapped by
$\widetilde{R}\ra\mathbb{P}^2$ to pencil of rational plane curves.

\begin{proposition}{\rm (See \cite[Proposition 3.8]{GSal})}
Let $\mathcal{B}$ be a generalized conic bundle over $\widetilde{R}$. Let $C$ be a section of $\mathcal{B}$. The pencil $\mathcal{B}$ induces a genus 1 fibration $\mathcal{E}_{\mathcal{B}}$ on the K3 surface $X$ which is the double cover of $\widetilde{R}$. The pull back of the curve $C$ is a section of the fibration $\mathcal{E}_{\mathcal{B}}$ if and only if $C$ is a branch curve if the double cover $X\ra \widetilde{R}$.
Moreover, all the elliptic fibrations on $X$ of type 1 are of this form.\end{proposition}

In the following subsections, for each value of $k$, we describe
the pencil of cubics used to construct the rational surface $R$
and then we summarize in tables the relations between the
generalized conic bundles on $R$ and the elliptic fibrations on
$X$. More precisely, by Proposition \ref{prop: when iota is
trivial on the basis}, each elliptic fibration on $X$ listed in
Table \ref{prop: possible fibrations} with only one exception for
each $k$ is induced by a generalized conic bundle on
$\widetilde{R}$. In the following tables, we associate to each of
these fibrations a generalized conic bundle inducing the elliptic
fibration. The elliptic fibration is described in the first and in
the last column: in the first it is given the number of the
fibration as given in Table of Proposition \ref{prop: possible
fibrations}, whereas in the last column, the reducible fibers of
the fibration and its Mordell--Weil group are described. In the
other columns, we provide a generalized conic bundle associated to
each elliptic fibration. The description of the generalized conic
bundle consists in giving the degree of the planes curves and the
list of the base points, since a generalized conic bundle is a
pencil of rational plane curves with some base points. We
distinguish between conic bundles (cb) and generalized conic
bundles which are not conic bundles (Gcb).

In order to associate to each generalized conic bundle an elliptic
fibration it suffices to find the reducible fibers of the
generalized conic bundle and to apply \cite[Theorem 5.3]{GSal}
which allows one to find the singular fibers of the elliptic
fibration associated to a conic bundle.

If some of the base points, say $p_1$ and $p_1'$, of the pencil of cubic curves that induces the rational elliptic surface are infinitely near, we say that a curve passes through $p_1$ and $p_1'$ to express that it shares the same tangent direction as the cubics in the rational elliptic pencil. This will be used in what follows.

From the following tables one obtains that it suffices to consider
generalized conic bundles which map to pencils of plane curves of
degree at most 3 in order to recover all the elliptic fibrations
listed in Proposition \ref{prop: possible fibrations}.

\subsection{The case $k=9$} The rational elliptic surface $\mathcal{E}$ is isomorphic to the blow-up of $\mathbb{P}^2$ in 9 points $p_1,p'_1,p''_1,p_2,p'_2, p'''_2, p_3,p'_3,p'''_3$, where $p_i,p'_i,p''_i$ are infinitely near points. Moreover the points $p_1,p_2,p_3$ are not collinear. We call $l_1$ the line through $p_1$ and $p_2$, the line joining $p_1$ and $p_3$ will be called $l_3$, while $l_2$ will be the line connecting $p_2$ and $p_3$. We assume that $l_i$ is tangent at $p_i$ to the cubics of the pencil that induces the elliptic fibration in $\mathcal{E}$.

$$\begin{array}{|c|c|c|c|c|}
\hline
n^o&\mbox{deg}&\mbox{ base points }&\mbox{ type }&\mbox{ elliptic fibrations }\\
\hline
9.3&1&p_1&\mbox{cb}&(III^*,I_6^*), MW=\Z/2\Z\\
\hline
9.4&2&p_1, p_1', p_1'', p_3&\mbox{cb}& (I_{12}^*, I_2),\ MW=\Z/2\Z\\
\hline 9.2&2&p_2,\mbox{ tg to }l_1, p_3\mbox{ tg to
}l_3&\mbox{Gcb}&(2II^*,I_2),\ MW=\{1\}\\ \hline
\end{array}$$

\subsection{The case $k=8$, $\delta=1$}
The rational elliptic surface $\mathcal{E}$ is isomorphic to the blow up of $\mathbb{P}^2$ in 9 points $p_1,p'_1, p^{''}_1, p_2, p'_2, p^{''}_2, p_3, p'_3, p_4$, such that $p_i$ and $p'_i$ and $p^{''}_i$, for $i=1,2,3$, are infinitely near points. These points are such that there is one triple of collinear points $l_3:=\{p_2,p_3, p_4\}.$ The line $l_2$ through $p_1,p_3$ is tangent at the point $p_1$, and the line $l_1$ through $p_1$ and $p_2$ is tangent at $p_2$ to all cubics of the pencil that gives the elliptic fibration on $\mathcal{E}$, and there is no other collinearity
relation.

$$\begin{array}{|c|c|c|c|c|}
\hline
\mbox{n}^o&\mbox{deg}&\mbox{ base points }&\mbox{ type }&\mbox{ elliptic fibrations }\\
\hline
8.3&1&p_1&\mbox{cb}&(III^*,I_4^*,I_2), MW=\Z/2\Z\\
\hline
8.5&1&p_2&\mbox{cb}& (I_6^*, I_2^*),\ MW=\Z/2\Z\\
\hline 8.4&2&p_1,\ p_3,\ p_3',\ p_4&\mbox{cb}&(I_{10}^*,2I_2),\
MW=\Z/2\Z\\
\hline 8.2&2&p_1,\mbox{ tg to }l_1,\ p_3,\ p_3' &\mbox{Gcb}&(II^*,III^*, I_2),\ MW=\{1\}\\
\hline
\end{array}$$

\subsection{The case $k=8$, $\delta=0$} The rational elliptic surface $\mathcal{E}$ is isomorphic to the blow up of $\mathbb{P}^2$ in 9 points $p_1,p'_1, p^{''}_1, p_2, p'_2, p^{''}_2, p_3, p'_3, p_4$, such that $p_i$ and $p'_i$
and $p^{''}_i$, for $i=1,2,3$, are infinitely near points. These
points are such that there is one triple of collinear points
$l_1:=\{p_1,p_2, p_4\}.$ The line $l_2$ through $p_1,p_3$ is
tangent at the point $p_1$, and the line $l_3$ through $p_2$ and
$p_3$ is tangent at $p_2$ to all cubics of the pencil that gives
the elliptic fibration on $\mathcal{E}$, and there is no other
collinearity relation.

$$\begin{array}{|c|c|c|c|c|}
\hline
\mbox{n}^o&\mbox{deg}&\mbox{ base points }&\mbox{ type }&\mbox{ elliptic fibrations }\\
\hline
8.6&1&p_1&\mbox{cb}&(2I_4^*),\ MW=\Z/2\Z\\
\hline
8.3&1&p_3&\mbox{cb}& (2III^*, 2I_2),\ MW=\Z/2\Z\\
\hline 8.5&1&p_4,&\mbox{cb}&(I_8^*, I_0^*),\
MW=\Z/2\Z\\
\hline
8.2&2&p_3,\mbox{ tg to }l_3,\ p_1,\ p_4&\mbox{Gcb}&(II^*, I_4^*),\ MW=\{1\}.\\
\hline 8.4&3&\begin{array}{l}p_1,\ p_1',\ p_1'',\ p_3\mbox{ tg to
}l_3,\\ \mbox{ with a node in }p_4\end{array}&\mbox{Gcb}&(I_{12}^*),\ MW=\{1\}.\\
\hline
\end{array}$$

\subsection{The case $k=7$}
The rational elliptic surface $\mathcal{E}$ is isomorphic to the
blow-up of $\mathbb{P}^2$ in 9 points
$p_1,p'_1,p''_1,p_2,p'_2,p_3,p'_3,p_4$ and $p_5$, where
$p_1,p'_1,p''_1$ are infinitely near points, and the same holds
for the pairs $p_2,p'_2$ and $p_3,p'_3$. These points are such
that there are two triples of collinear points
$l_2:=\{p_1,p_4,p_3\}$ and $l_3:=\{p_2,p_3,p_5\}$ and no other
collinearity relation. We call $l_1$ the line through $p_1$ and
$p_2$, and we assume that it is tangent to the cubics at $p_1$.

$$\begin{array}{|c|c|c|c|c|}
\hline
\mbox{n}^o&\mbox{deg}&\mbox{ base points }&\mbox{ type }&\mbox{ elliptic fibrations }\\
\hline
7.9&1&p_1&\mbox{cb}&(I_4^*,I_2^*,I_2),\ MW=\Z/2\Z\\
\hline
7.3&1&p_2&\mbox{cb}& (III^*,I_2^*,2I_2),\ MW=\Z/2\Z\\
\hline 7.8&1&p_4,&\mbox{cb}&(I_6^*, I_0^*,I_2),\
MW=\Z/2\Z\\
\hline
7.7&2&p_2,\ p'_2,\ p_1,\ p_5&\mbox{cb}&(I_8^*,3I_2),\ MW=\Z/2\Z.\\
\hline
7.5&2&p_1,\ p'_1,\ p_3\mbox{ tg to }l_3&\mbox{Gcb}&(2III^*,I_2),\ MW=\{1\}.\\
\hline
7.2&2&p_2\mbox{ tg to }l_1,\ p_3,\ p_4&\mbox{Gcb}&(II^*, I_2^*, I_2),\ MW=\{1\}.\\
\hline
7.4&2&p_2\mbox{ tg to }l_1,\ p_4,\ p_5&\mbox{Gcb}&(III^*,I_4^*),\ MW=\{1\}.\\
\hline
7.6&3&\begin{array}{l}p_2\mbox{ tg to }l_1,\ p_3,\ p_3',\ p_5,\ \\
\mbox{ with a node in }p_4\end{array}&\mbox{Gcb}&(I_{10}^*, I_2),\ MW=\{1\}.\\
\hline
\end{array}$$

\subsection{The case $k=6$} The rational elliptic surface $\mathcal{E}$ is isomorphic to the blow up of $\mathbb{P}^2$ in 9 points $p_1,p'_1, p_2, p'_2, p_3, p'_3, p_4, p_5, p_6$, such that $p_i$ and $p'_i$, for $i=1,2,3$, are infinitely near points. These points are such that there are three triple of collinear points
$l_1:=\{p_1,p_2, p_4\}, l_2:=\{p_1,p_5,p_3\},
l_3:=\{p_2,p_3,p_6\}$, and no other collinearity relation (see
\cite{Pastro}).

$$\begin{array}{|c|c|c|c|c|}
\hline
\mbox{n}^o&\mbox{deg}&\mbox{ base points }&\mbox{ type }&\mbox{ elliptic fibrations }\\
\hline
6.10&1&p_1&\mbox{cb}&(2I_2^*,2I_2),\ MW=\Z/2\Z\\
\hline
6.9&1&p_4&\mbox{cb}&(I_4^*,I_0^*,2I_2),\
MW=\Z/2\Z\\
\hline
6.4&2&p_1,\ p_2,\ p_5,\ p_6&\mbox{cb}&(III^*, I_0^*,3I_2),\ MW=\Z/2\Z\\
\hline
6.7&2&p_1,\ p'_1,\ p_2,\ p_5&\mbox{cb}&(I_6^*,4I_2),\ MW=\Z/2\Z\\
\hline
6.8&2&p_1,\ p_1',\ p_5,\ p_6&\mbox{cb}&(I_4^*,I_2^*),\ MW=\{1\}\\
\hline
6.2&3&\begin{array}{l}p_1,\ p_1',\ p_3\mbox{ tg to }l_3,\ p_6,\ \\
\mbox{ with a node in }p_4\end{array}&\mbox{Gcb}&(II^*, I_0^*, 2I_2),\ MW=\{1\}\\
\hline
6.6&3&\begin{array}{l}p_1,\ p_1',\ p_3\mbox{ tg to }l_2,\ p_6,\ \\
\mbox{ with a node in }p_4\end{array}&\mbox{Gcb}&(I_6^*,I_0^*),\ MW=\{1\}\\
\hline
6.5&3&\begin{array}{l}p_1,\ p_1',\ p_3\mbox{ tg to }l_3,\ p_5,\ \\
\mbox{ with a node in }p_4\end{array}&\mbox{Gcb}&(I_8^*,2I_2),\ MW=\{1\}\\
\hline
6.3&3&\begin{array}{l}p_1,\ p_3\mbox{ tg to }l_3,\ p_5,\ p_6\\
\mbox{ with a node in }p_4\end{array}&\mbox{Gcb}&(III^*,I_2^*, I_2),\ MW=\{1\}\\
\hline
\end{array}$$

\subsection{The case $k=5$}  The rational elliptic surface $\mathcal{E}$ is isomorphic to the blow up of $\mathbb{P}^2$ in 9 points $p_1,p'_1, p_2, p'_2, p_3,p _4,p_5,p_6, p_7$, such that $p_1$ and $p'_1$ are infinitely near points, as are $p_2$ and $p'_2$. These points are such that there are three triple of collinear points
$l_1:=\{p_1,p_2, p_3\}, l_2:=\{p_1,p_4,p_5\}, l_3:=
\{p_2,p_6,p_7\}$, and no other collinearity relation (see
\cite{Salgado}). We call $q$ the intersection $l_2\cap l_3$.

$$\begin{array}{|c|c|c|c|c|}
\hline
\mbox{n}^o&\mbox{deg}&\mbox{ base points }&\mbox{ type }&\mbox{ elliptic fibrations }\\
\hline
5.9&1&p_1&\mbox{cb}&(I_2^*, I_0^*, 3I_2),\ MW=\Z/2\Z\\
\hline
5.7&1&p_3&\mbox{cb}&(I_4^*,5I_2),\
MW=\Z/2\Z\\
\hline 5.8&1&q&\mbox{Gcb}&(2I_2^*,I_2),\
MW=\{1\}\\
\hline
5.4&2&p_2,\ p_3,\ p_4,\ p_5&\mbox{cb}&(III^*,6I_2),\ MW=\Z/2\Z\\
\hline
5.3&2&p_1,\ p_3,\ p_6,\ q&\mbox{Gcb}&(III^*, I_0^*,2I_2),\ MW=\{1\}\\
\hline 5.6&2&p_1,\ p_1',\ p_2,\ q&\mbox{Gcb}&(I_4^*,I_0^*, I_2),\ MW=\{1\}\\
\hline
5.5&2&p_1,\ p_1',\ p_3,\ q&\mbox{Gcb}&(I_6^*,3I_2),\ MW=\{1\}\\
\hline
5.2&3&\begin{array}{l}p_1,\ p_1',\ p_6,\ p_7,\ q\\
\mbox{ with a node in }p_3\end{array}&\mbox{Gcb}&(II^*,5I_2),\ MW=\{1\}\\
\hline
\end{array}$$

\subsection{The case $k=4$, $\delta=0$}
The rational elliptic surface $\mathcal{E}$ is isomorphic to the
blow up of $\mathbb{P}^2$ in 9 points $p_1,..., p_9$. These points
are such that there are three triple of collinear points
$l_1:=\{p_1,p_2, p_3\}$, $l_2:=\{p_4,p_5,p_6\}$,
$l_3:=\{p_7,p_8,p_9\}$, and no other collinearity relation (see
\cite{Fusi}). The lines $l_1$, $l_2$ and $l_3$ meet in a
unique point $q$.

$$\begin{array}{|c|c|c|c|c|}
\hline
\mbox{n}^o&\mbox{deg}&\mbox{ base points }&\mbox{ type }&\mbox{ elliptic fibrations }\\
\hline
4.2&1&p_1&\mbox{cb}&(I_2^*,6I_2),\ MW=\Z/2\Z\\
\hline 4.3&1&q&\mbox{Gcb}&(3I_0^*),\
MW=\{1\}\\
\hline
\end{array}$$

\subsection{The case $k=4$, $\delta=1$}
The rational elliptic surface $\mathcal{E}$ is isomorphic to the
blow up of $\mathbb{P}^2$ in 9 points $p_1,...,p_8, p'_8$, such
that $p_8$ and $p'_8$ are infinitely near points. These points are
such that there are three triple of collinear points
$l_1:=\{p_1,p_2, p_3\}$, $l_2:=\{p_4,p_5,p_8\}$,
$l_3:=\{p_6,p_7,p_8\}$, and no other collinearity relation (see
\cite{Fusi}). We call $q_1$, resp. $q_2$ the intersection point of
the first line $l_1$ with the the line $l_2$, resp. $l_3$.

$$\begin{array}{|c|c|c|c|c|}
\hline
\mbox{n}^o&\mbox{deg}&\mbox{ base points }&\mbox{ type }&\mbox{ elliptic fibrations }\\
\hline
4.6&1&p_1&\mbox{cb}&(2I_0^*,4I_2),\ MW=\Z/2\Z\\
\hline
4.5&1&p_4&\mbox{cb}&(I_2^*,6I_2),\
MW=\Z/2\Z\\
\hline 4.4&1&q_1&\mbox{Gcb}&(I_2^*,I_0^*,2I_2),\
MW=\{1\}\\
\hline 4.3&2&q_1,\ p_1,\ p_4,\ p_5&\mbox{Gcb}&(I_4^*,4I_2),\ MW=\{1\}\\
\hline
4.2&2&q_1,\ p_4,\ p_6,\ p_7&\mbox{Gcb}&(III^*,5I_2),\ MW=\{1\}\\
\hline
\end{array}$$

\subsection{The case $k=3$} The rational elliptic surface $\mathcal{E}$ is isomorphic to the blow up of $\mathbb{P}^2$ in 9 distinct points $p_1,...,p_9$, such that there are three triple of collinear points
$l_1:=\{p_1,p_2, p_3\}$, $l_2:=\{p_4,p_5,p_6\}$,
$l_3:=\{p_7,p_8,p_9\}$, and no other collinearity relation (see
\cite{Fusi}). We call $q_1,q_2$ and $q_3$ the three intersection
points of the three pairs of these lines, with the assumption that
$q_3=l_2\cap l_3$.

$$\begin{array}{|c|c|c|c|c|}
\hline
\mbox{n}^o&\mbox{deg}&\mbox{ base points }&\mbox{ type }&\mbox{ elliptic fibrations }\\
\hline
3.4&1&p_1&\mbox{cb}&(I_0^*,7I_2),\ MW=\Z/2\Z\\
\hline 3.3&1&q_1&\mbox{Gcb}&(2I_0^*,3 I_2),\
MW=\{1\}\\
\hline 3.2&2&p_1,\ p_2,\ p_4,\ q_3&\mbox{Gcb}&(I_2^*,5I_2),\ MW=\{1\}\\
\hline
\end{array}$$

\subsection{The case $k=2$} The rational elliptic surface $\mathcal{E}$ is isomorphic to the blow up of $\mathbb{P}^2$ in 9 points, such that 3 of them lie on a line $l$, say $p_1,p_2,p_3$, and the remaining 6, namely, $p_4,...,p_9,$ lie on a conic $Q$ (see \cite{Fusi}). We denote by $q_i$, for $i=1,2$, the two intersection points $Q\cap l$.

$$\begin{array}{|c|c|c|c|c|}
\hline
\mbox{n}^o&\mbox{deg}&\mbox{ base points }&\mbox{ type }&\mbox{ elliptic fibrations }\\
\hline
2.3&1&p_4&\mbox{cb}&(10 I_2),\ MW=\Z/2\Z\\
\hline
2.2&1&q_1&\mbox{Gcb}&(I_0^*,6I_2),\ MW=\{1\}\\
\hline
\end{array}$$

\subsection{The case $k=1$}
The rational elliptic surface $\mathcal{R}$ is isomorphic to the
blow up of $\mathbb{P}^2$ in 9 points $p_1,\ldots, p_9$ in general
position. There are twelve nodal cubics in the pencil of cubics
through $p_1,\ldots, p_9$. We choose one of them as branching
fiber of the double cover and we denote by $q$ its singular point.

$$\begin{array}{|c|c|c|c|c|}
\hline
\mbox{n}^o&\mbox{deg}&\mbox{ base points }&\mbox{ type }&\mbox{ elliptic fibrations }\\
\hline
1.2&1&q&\mbox{Gcb}&(9 I_2),\ MW=\{1\}\\
\hline
\end{array}$$

\section{Equations}\label{sec: equations}
In the paper \cite{wine2}, the authors give an algorithm to compute the
Weierstrass equations of certain elliptic fibrations on K3
surfaces which are double cover of rational elliptic surfaces. The aim of this section is to recall such an algorithm, and
to compute explicitly some of these equations. We observe that for
each elliptic fibration described in the previous section one can
apply the algorithm and therefore find explicitly the Weierstrass
equations.

\subsection{An example}
Let us consider the rational elliptic surface $\mathcal{R}_k$, for $k=1,\ldots, 9$,
associated to the pencil of cubics $\mathcal{P}_k$ as in
Proposition \ref{prop: equations of RES}. Let $X$ be the K3
surface branched on two fibers of the rational elliptic surface
$R_k$. We chose as branch fibers the one over $\infty$ and another
one, over $\mu_1$. So the K3 surface is obtained as double cover
of $\mathbb{P}^2$ branched on the sextic $\mathcal{S}_k$
\begin{multline} (xyz+a_3z^3+a_2y^3)(x^2y+y^2z+z^2x+a_9xy^2+a_8x^2z+a_7yz^2+a_6z^3+a_5y^3+a_4x^3+\\ \mu_1(xyz+a_3z^3+a_2y^3)),\end{multline}
with $a_2=\ldots=a_k=0$.

If $k\leq 3$, the point $(1:0:0)$ is the singular point of the
fiber over $\infty$ and it is not a base point of the pencil
$\mathcal{P}_k$. If $k>3$, then $(1:0:0)$ is a base point of the
pencil $\mathcal{P}_k$. For every $k$, the pencil of lines through
$(1:0:0)$ induces a generalized conic bundle on $\widetilde{R}$,
which is also a conic bundle on $R$ if $k>3$.

In order to find the Weierstrass equation of the elliptic
fibration induced by this generalized conic bundle, one has to
consider the pencil of lines through $(1:0:0)$, $y=mz$. Then one
intersects it with the branch sextic $\mathcal{S}_k$ obtaining the
following equation for the elliptic fibration on the K3 surface
$X$:
\begin{eqnarray*}\begin{array}{l}Y^2=z^2\left(m+z(a_3+a_2m^3)\right)\left(a_4+z(m+a_8)+\right.\\ \left.z^2(1+a_9m^2+\mu_1m)+
z^3\left(m^3(a_5+\mu_1a_2)+
m^2+a_7m+a_6+\mu_1a_3\right)\right).\end{array}\end{eqnarray*}

By the change of coordinates $Y\mapsto Yz$, one obtains $Y^2$
equals to a polynomial of degree $4$ in $z$ with a section (if $a_2=a_3=0$ the section is at infinity otherwise it is $(z(m),Y(m))=(-m/(a_3+a_2 m^3),0)$. So we obtain the equation of
an elliptic fibration. If $a_2=a_3=0$, after the change of
coordinates
\begin{eqnarray*}\begin{array}{lll}Y&\mapsto&
Y/\left(m\left(m^3(a_5+\mu_1a_2)+ m^2+a_7m+a_6+\mu_1a_3\right)^2\right),\\
z&\mapsto &z/\left(m\left(m^3(a_5+\mu_1a_2)+
m^2+a_7m+a_6+\mu_1a_3\right)\right)\end{array}\end{eqnarray*} one
obtains the Weierstrass form. If $k>3$, putting $a_2=\ldots=a_k=0$
one obtains a Weierstrass equation for the fibrations $n^o$
$(k.h)$ in Proposition \ref{prop: possible fibrations},  for the
following pairs of values
$(k,h)\in\{(4,6),(5,9),(6,10),(7,9),(8,3),(9,3)\}$. For $k=1$
(resp. $k=2$, $k=3$) this is an equation for the fibration 1.2
(resp. 2.2, 3.3) in Proposition \ref{prop: possible fibrations}.

\subsection{The algorithm}\label{subsection: the algorithm} The aim of this algorithm is to generalize the previous
computation in order to be able to write the Weierstrass equation
of all the elliptic fibration classified in Proposition \ref{prop:
possible fibrations}

\noindent\textbf{Setup.} Let $V$ be a K3 surface obtained by a
base change of order 2 from a rational elliptic surface $R$.
Therefore, $V$ can be described as double cover of $\mathbb{P}^2$
branched on the union of two (possibly reducible) plane cubics
from the pencil determining $R$. It has an equation of the form
\begin{equation*}w^2=f_3(x_0:x_1:x_2)g_3(x_0:x_1:x_2).\end{equation*}
Let $\mathcal{B}$ be a (generalized) conic bundle on
$\widetilde{R}$ whose curves are parametrized by $\tau$. Pushing
forward to $\mathbb{P}^2$, $\mathcal{B}$ is given by a pencil of
plane rational curves with equation $h(x_0:x_1:x_2,\tau)$.  The
polynomial $h(x_0:x_1:x_2,\tau)$ is homogeneous in $x_0$, $x_1$,
$x_2$, say of degree $e\geq 1$ and linear in $\tau$.

The adjunction formula implies that every curve with equation
$h(x_0:x_1:x_2,\tau)$ meets both of the branch curves (the proper
transforms on $\widetilde{R}$ of) $f_3=0$ and $g_3=0$ in two
additional points not blown up by $\widetilde{R}\ra\mathbb{P}^2$.
It therefore meets (the proper transform of) their union
$f_3g_3=0$ in four points not blown up by
$\widetilde{R}\ra\mathbb{P}^2$.  So the preimage in $V$ is the
double cover of a rational curve branched over $4$ points, i.e.,
the standard presentation of an elliptic curve.  For general
$\tau$, we must find an isomorphism of the curve
$h(x_0:x_1:x_2,\tau) =0$ with $\mathbb{P}^1$, and extract the
images of the four intersection points with $f_3g_3 =0$.

When all curves in the conic bundle have a basepoint of degree
$e-1$, the projection away from this point provides the required
isomorphism of the curve $h(x_0:x_1:x_2,\tau) =0$ with
$\mathbb{P}^1$. Up to acting by $\PGL_3(\C)$, we may assume that
this point is $(0:1:0)$.

\noindent\textbf{Algorithm.}
\begin{enumerate}
\item Compute the resultant of the polynomials
$f_3(x_0:x_1:x_2)g_3(x_0:x_1:x_2)$ and $h(x_0:x_1:x_2,\tau)$ with
respect to the variable $x_1$. The result is a polynomial
$r(x_0:x_2,\tau)$ which is homogeneous in $x_0$ and $x_2$,
corresponding to the images of \textit{all} of the intersection
points $\{f_3g_3=0 \} \cap \{h_\tau=0\}$ after projection from
$(0:1:0)$.

\item Since $\mathcal{B}$ is a conic bundle, $r(x_0:x_2,\tau)$
will be of the form $a(x_0:x_2,\tau)^2b(x_0:x_2,\tau)c(\tau)$,
where $a$ and $b$ are homogeneous in $x_0$ and $x_2$, the degree
of $a$ depends upon $e$ and the degree of $b$ in $x_0$ and $x_2$
is 4. \item The equation of $V$ is now given by
$w^2=r(x_0:x_2,\tau)$, which is birationally equivalent to
\begin{equation}\label{eq: genus 1 fibration by conic bundles}w^2=c(\tau)b(x_0:x_2,\tau),\end{equation} by the change of coordinates
$w\mapsto wa(x_0:x_2,\tau)$. Since for almost every $\tau$, the
equation \eqref{eq: genus 1 fibration by conic bundles} is the
equation of a $2:1$ cover of $\mathbb{P}^1_{(x_0:x_2)}$ branched
in 4 points, \eqref{eq: genus 1 fibration by conic bundles} is the
equation of the genus 1 fibration on the K3 surface $V$ induced by
the conic bundle $\mathcal{B}$. \item If there is a section of
fibration \eqref{eq: genus 1 fibration by conic bundles}, then it
is possible to obtain the Weierstrass form by standard
transformations.\end{enumerate}

\begin{rem}{\rm There are several conic bundles whose general member
cannot be parametrized by lines. An algorithm for some of them is
described in \cite{wine2}, but here we do not need it, since all
the conic bundles listed in Section \ref{sec: geometric interpretations} are at most of degree
3.}\end{rem}

\section{Specializations}\label{sec: specializations}

\subsection{Specialization of a 1-dimensional family of K3 surfaces}

Among the rational elliptic surfaces listed in Section \ref{sec:
RES} the one  with the smallest Mordell--Weil group is the
extremal rational elliptic surface with a fiber of type $I_9$, 3
fibers of type $I_1$ and $MW=\Z/3\Z$.

As already noticed, it is induced by a pencil of cubics
$\mathcal{P}_9$ on $\mathbb{P}^2$.
By standard transformations one obtains the Weierstrass equation
of the rational elliptic surface:
\begin{equation}\label{eq: Weierstrass RES I9}Y^2=X^3+X\left(-\frac{1}{2^43}\mu^4-\frac{\mu}{2}\right)+\frac{1}{2^53^3}\mu^6+\frac{\mu^3}{2^33}+\frac{1}{2^2}.\end{equation}
The discriminant is $\frac{1}{16}\left(3+\mu\right)\left(\mu^2-3\mu+9\right)$.\\

Now let us consider the K3 surface $X_{\mu_1}$ obtained by a base
change of order 2 of this rational elliptic surface branched over the fiber of type
$I_9$ (corresponding to $\mu=\infty$) and a generic smooth fiber,
say the fiber corresponding to $\mu_1$.

Its Weierstrass equation is obtained by \eqref{eq: Weierstrass RES
I9} by substituting $\mu$ with $\tau^2+\mu_1$. For a
generic choice of $\mu_1$ this corresponds to an elliptic
fibration with $I_{18}+6I_1$ as reducible fibers and Mordell
Weil group given by $\Z/3\Z$. The transcendental lattice of K3 surfaces
in this family is $U\oplus\langle 2\rangle$, indeed they are the
K3 surfaces with one involution acting trivially on the
N\'eron--Severi group and fixing 9 rational curves and 1 curve of
genus 1.

These K3 surfaces specialize to several K3 surfaces whose
transcendental lattice has rank 2 (i.e., whose Picard number is
20). For example, since $\langle 2d\rangle$ is primitively
embedded in $U$, the K3 surfaces whose transcendental lattice is
$\langle 2d\rangle\oplus\langle 2 \rangle$, $d>0$, are special
members of the same family.

Let us consider one of this specializations in details. Consider the plane conic $\mathcal{C}:=\{ xz+y^2\}$. It intersects
the cubic of the pencil corresponding to $\overline{\mu}$ in
$(0:0:1)$ with multiplicity 3, in $(1:0:0)$ with multiplicity 1
and in $(-\overline{\mu}:\pm\sqrt{\overline{\mu}}:1)$ with
multiplicity 1. This conic is a bisection of the rational elliptic
fibration. Indeed, it intersects the generic cubic of the pencil in
exactly 2 points which are not base points. If $\overline{\mu}=0$,
then the points $(0:0:1)$ and
$(-\overline{\mu}:\pm\sqrt{\overline{\mu}}:1)$ collapse to the
same point, thus the conic $\mathcal{C}$ intersects the cubic of
the pencil corresponding to $\mu=0$ in $(0:0:1)$ with multiplicity
5 and $(1:0:0)$ with multiplicity 1.

If now one considers the K3 surfaces $X_{\mu_1}$ obtained by the
base change branched on $\mu=\infty$ and $\mu_1$, generically the
bisection of the rational elliptic surface induced by $\mathcal{C}$ induces a bisection
of the elliptic fibration on the K3 surface. This does not happen if $\mu_1$ is 0. Indeed,
in this case, the bisection of the rational surface splits in the double cover
so that it induces, on the elliptic fibration on $X_0$, two distinct
sections. This is due to the fact that $X_0$ is the double cover
of $\mathbb{P}^2$ branched on the sextic $xyz(x^2y+y^2z+z^2x)$
which intersects the conic $\mathcal{C}$ in $(0:0:1)$ with
multiplicity 8 and $(0:0:1)$ with multiplicity 4, in particular
always with an even multiplicity, so it splits in the double
cover.

From the above discussion, we have that if the base change from the rational surface to the
K3 surface is branched on $\mu=0$ and $\mu=\infty$, then the Picard
number of the K3 surface $X_0$ is not 19, but 20. So all the
fibrations on the K3 surface $X_0$ specialize. This suggests a
method to determine some specializations: we have many elliptic
fibrations on the same K3 surface. By requiring that the
discriminant of one of these fibrations has zeros with higher
multiplicity, one obtains values of $\mu_1$ for which the K3
surface specializes, and thus its N\'eron--Severi group is larger than
the expected one. One is now able to find different
specializations considering different elliptic fibrations. In the
following we apply this idea to different families of K3
surfaces.

\subsubsection{The elliptic fibrations on $X_{\mu_1}$ and their specilizations}

For a generic choice of $\mu_1$, there are four types of elliptic
fibrations on $X_{\mu_1}$: one comes from the rational elliptic
surface, the other from generalized conics bundles on
$\widetilde{R}$.
\begin{enumerate}\item The elliptic fibration $\mathcal{E}_1$ coming from the rational elliptic surface has $I_{18}+6I_1$ as singular fibers and $MW=\Z/3\Z$. Its equation is
$$
Y^2=X^3+\left(\frac{1}{48}(\tau^2+\mu_1)(\tau^6+3\tau^4\mu_1+3\tau^2\mu_1^2+\mu_1^3+24)\right)X-\frac{1}{864}(\tau^2-\mu_1)^{6}+\frac{1}{4}+\frac{1}{24}(\tau^2-\mu_1)^3.
$$
and the discriminant is
$$\Delta_1:= \frac{1}{16}\left(\tau^2+\mu_1+3\right)\left(\mu_1^2+2\tau^2\mu_1-3\mu_1+\tau^4-3\tau^2+9\right);$$
\item The elliptic fibration $\mathcal{E}_2$ is associated to the
conic bundle $y=mx$ (lines through $(0:0:1)$), its singular fibres
are $III^*+I_{6}^*+3I_1$, $MW=\Z/2\Z$ and its equation is $$
Y^2=X^3-\frac{1}{3}m^3(m^3+2m^2\mu_1+\mu_1^2m-3)X-\frac{1}{27}m^5(m+\mu_1)(2m^3+4m^2\mu_1+2\mu_1^2m-9)
$$
$$\Delta_2:=-m^9(-4+m^3+2m^2\mu_1+\mu_1^2m);$$
\item The elliptic fibration $\mathcal{E}_3$ is associated to the
conic bundle $y^2+bxy+xz$, its singular fibres are
$I_{12}^*+I_2+4I_1$, $MW=\Z/2\Z$ and its equation is
$$Y^2=X^3+A_3X+B_3$$
$$A_3:=-(1/3)b^6+(4/3)b^3+(2/3)b^5\mu_1-1/3-(4/3)\mu_1 b^2-(1/3)\mu_1^2b^4$$
$$B_3:=(1/27)(-\mu_1 b^2-2+b^3)(2b^6-4b^5\mu_1+2\mu_1^2b^4-8b^3+8\mu_1 b^2-1)$$
$$\Delta_3:=-b^2(-\mu_1 b^2-4+b^3)(-\mu_1+b);$$
\item The elliptic fibration $\mathcal{E}_4$ is associated to the
generalized conic bundle $ax^2+yz$, its singular fibres are
$2II^*+I_{2}+2I_1$, $MW=\{1\}$ and its equation is
$$Y^2=X^3-\frac{1}{3}\mu_1^2a^4X+\frac{a^5}{27}(27a^2-54a-2\mu_1^3a+27)$$
$$\Delta_4:=a^{10}(-1+a)^2(27a^2-4\mu_1^3a-54a+27).$$
\end{enumerate}
A very natural specialization for elliptic fibrations is obtained by
requiring that certain singular fibers collapse to a unique one.

By $\mathcal{E}_1$ one obtains that possible specializations,
under which the Picard number jumps to $20$ are obtained by
requiring $$\mu_1\in\left\{-3,\frac{3-3\sqrt{3}i}{2},
\frac{3+3\sqrt{3}i}{2}\right\}$$ (in this case the second branch
fiber is a fiber of type $I_1$ on the rational elliptic surface
and thus gives a fiber of type $I_2$ on the K3 surface).

By $\mathcal{E}_2$ one obtains that the Picard numbers jump for
same values of $\mu_1$.

By $\mathcal{E}_4$ and $\mathcal{E}_3$ one obtains that the Picard
numbers jump for $\mu_1=0$ and the same values found for
$\mathcal{E}_1$.
\subsubsection{Considering explicitly the specializations}
If $$\mu_1\in\left\{-3,\frac{3-3\sqrt{3}i}{2},
\frac{3+3\sqrt{3}i}{2}\right\},$$ one obtains that $\mathcal{E}_1$
has a new reducible fiber of type $I_2$ (obtained gluing together
two fibers of type $I_1$). Thus the trivial lattice of the
specialized elliptic fibration is $U\oplus A_{17}\oplus A_1$. The
3-torsion section generating the Mordell Weil group does not
change, but there cannot be other torsion sections (by
\cite{Shim} the torsion part of the Mordell Weil group of an
elliptic fibration with trivial lattice $U\oplus A_{17}\oplus A_1$
is either $\{1\}$ or $\Z/3\Z$). Moreover, by the Shioda-Tate formula,
there is no  section of infinite order, since the rank of the
trivial lattice is the maximum admitted. Hence we have a set of
generators for the N\'eron--Severi of the specialized surface, and
thus we can also compute the transcendental lattice: it is
$\langle 2\rangle\oplus\langle 2\rangle$.

Let us consider the specialization of the elliptic fibrations
$\mathcal{E}_2$, $\mathcal{E}_3$, $\mathcal{E}_4$ if
$\mu_1\in\left\{-3,\frac{3-3\sqrt{3}i}{2},
\frac{3+3\sqrt{3}i}{2}\right\}.$ In all the cases an extra
reducible fiber appears. Indeed, the singular fibers of the
fibration $\mathcal{E}_2$ become $III^*+I_{6}^*+I_2+I_1$, the ones
of $\mathcal{E}_3$ become $I_{12}^*+2I_2+2I_1$ and
$\mathcal{E}_4$ become $2II^*+2I_2$.

The unique other interesting value of $\mu_1$ found before is
$\mu_1=0$. Let us denote by $X_0$ the K3 surface obtained for
$\mu_1=0$.

The elliptic fibration $\mathcal{E}_4$ has a new reducible fiber
if $\mu_1=0$. The reducible fibers are $2II^*+IV$ and thus the
trivial lattice is $U\oplus E_8\oplus E_8\oplus A_2$. The
Mordell--Weil is trivial and so $NS(X_0)\simeq U\oplus E_8\oplus
E_8\oplus A_2$ and $T_{X_0}\simeq A_2$. If $\mu_1=0$, the equation
of $\mathcal{E}_4$ is $$Y^2=X^3+\frac{a^5}{27}(27a^2-54a+27)$$
which is clearly an isotrivial fibration whose generic fiber is
isomorphic to the elliptic curve with complex multiplication of
order 3. Let us denote by $\alpha$ the order 3 automorphism
induced on $X_0$ by the complex multiplication on the fibers; it
acts trivially on the basis of the fibration and with order 3 on
each fiber.

We observe that if $\mu_1=0$, the elliptic fibration
$\mathcal{E}_2$ has the equation
$$
Y^2=X^3-\frac{1}{3}m^3(m^3-3)X-\frac{1}{27}m^6(2m^3-9)
$$
which admits the new automorphism of order 3, $(x,y,m)\mapsto
(x,y,\zeta_3 m)$. It coincides with $\alpha$. The singular fibers
of this fibration are unchanged, but now one can view this
fibration as $3:1$ cover of the rational elliptic fibration whose
equation is
$$
Y^2=X^3-\frac{1}{3}n(n-3)X+-\frac{1}{27}n^2(2n-9)
$$
and whose singular fibers are $I_2^*+III+I_1$ and $MW=\Z$. The
branch fiber of the triple cover are $I_2^*$ and $III$. The
fibration $\mathcal{E}_2$ over $X_0$, admits an extra section of
infinite order, induced by the one of the rational elliptic
surface .

The action of $\alpha$ on the  fibration $\mathcal{E}_3$ is
similar to the previous one, $\alpha$ acts on the basis of the
fibration. But, moreover, one observes that the discriminant
$\Delta_3$ changes if $\mu_1=0$, and indeed it is $-b^3(-4+b^3).$
So the trivial lattice of $\mathcal{E}_3$ specialized to $X_0$ is
$U\oplus D_{16}\oplus A_2.$ The Mordell Weil group remains
unchanged.

The results obtained are summarized in the following Table.
$$
\begin{array}{|c|c|c|}
\hline
\mu_1&T_{X_{\mu_1}}&\mbox{ fibrations }\\
\hline
\{-3,\frac{3-3\sqrt{3}i}{2}, \frac{3+3\sqrt{3}i}{2}\}&\langle2\rangle^2 &\begin{array}{|c|c|c|}\mathcal{E}_i&\mbox{ fibers }&MW\\
\hline
\mathcal{E}_1&I_{18}+I_2+4I_1&\Z/3\Z\\
\mathcal{E}_2&III^*+I_{6}^*+I_2+I_1
&\Z/2\Z\\
\mathcal{E}_3&I_{12}^*+2I_2+2I_1&\Z/2\Z\\
\mathcal{E}_4&2II^*+2I_2&\{1\}
\end{array}\\
\hline 0&A_2(-1)&
\begin{array}{|c|c|c|}\mathcal{E}_i&\mbox{ fibers }&MW\\
\hline
\mathcal{E}_1&I_{18}+6I_1&\Z/3\Z\times \Z\\
\mathcal{E}_2&III^*+I_6^*+3I_1&\Z/2\Z\times \Z\\
\mathcal{E}_3&I_{12}^*+I_3+3I_1&\Z/2\Z\\\
\mathcal{E}_4&2II^*+IV&\{1\}\\

\end{array}\\
\hline
\end{array}
$$

In particular we obtained the``two most algebraic K3 surfaces" by
Vinberg, see \cite{V}, as specializations of our K3 surface. We
recall that the elliptic fibrations on a rigid K3 surface whose
transcendental lattice is either $\langle
2\rangle\oplus\langle2\rangle$ or $A_2(-1)$ were classified by
Nishiyama, see \cite{Nis}.

\subsection{Specializations of 2-dimensional families}
\subsubsection{$\delta=1$}
We consider the 1-dimensional family of rational elliptic surfaces, parametrized by $a$, associated to the pencil of cubics:
$$x^2y+y^2z+z^2x+axy^2+\mu xyz.$$
The Weierstrass equation of the rational elliptic fibration
$\mathcal{E}_R:R_a\ra\mathbb{P}^1_{\mu}$ is
\begin{equation}\label{eq: RES delta=1}Y^2=X^3+\left(-\frac{1}{48}\mu^4+\frac{1}{6}\mu^2a-\frac{1}{3}a^2-\frac{1}{2}\mu\right)X-\frac{1}{864}\mu^6+\frac{1}{72}\mu^4a-\frac{1}{24}\mu^3-\frac{1}{18}\mu^2a^2+\frac{1}{6}\mu
a+\frac{2}{27}a^3-\frac{1}{4}\end{equation} and the discriminant
is
$$\Delta:=-\frac{1}{16}\mu^4a+\frac{1}{16}\mu^3+\frac{1}{2}\mu^2a^2-\frac{9}{4}\mu
a-a^3+\frac{27}{16}.$$

To obtain the K3 surface $X_{a,\mu_1}$ with the elliptic fibration
$\mathcal{E}_X:X_{a,\mu_1}\ra\mathbb{P}^1_{\tau}$ we apply the
base change $\mu=\tau^2+\mu_1$ which is branched at infinity and
at $\mu=\mu_1$.

On $X_{a,\mu_1}$ we have 5 elliptic fibrations: one of them is
induced by $\mathcal{E}_R$ after the base change and its
Weierstrass equation is obtained directly by \eqref{eq: RES
delta=1}, the others are associated to generalized conic bundles
and their equations can be find applying the method described in
Section \ref{subsection: the algorithm}. Since in the following we
are mainly interested in the discriminant of these fibrations, we
only write the equations of the conic bundles and the discriminant
of the associated elliptic fibrations.

The 5 elliptic fibrations are:\begin{enumerate}\item
$\mathcal{E}_1$ induced by the elliptic fibration $\mathcal{E}_R$
on $R_a$, whose discriminant is
\begin{eqnarray*}\label{eq: E1 2-dim d1}
\begin{array}{l}
\Delta_1:=-\frac{1}{16}(\tau^2+\mu_1)^4a+\frac{1}{16}(\tau^2+\mu_1)^3+\frac{1}{2}(\tau^2+\mu_1)^2a^2-\frac{9}{4}(\tau^2+\mu_1)
a-a^3+\frac{27}{16}
\end{array}
\end{eqnarray*}
the singular fibers are $I_{16}+8I_1$ and $MW=\Z$; \item
$\mathcal{E}_2$ induced by the conic bundle $x-my$, whose
discriminant is
\begin{eqnarray*}\label{eq: E2 2-dim d1}\begin{array}{l}
\Delta_2:=m^{10}(a+m)^2\left(4m^2a-1+4m^3-m^2\mu_1^2-2m\mu_1\right)
\end{array}\end{eqnarray*}
the singular fibers are $III^*+I_4^*+I_2+3I_1$ and $MW=\Z/2\Z$;
\item $\mathcal{E}_3$ induced by the conic bundle $y-mz$, whose
discriminant is
\begin{eqnarray*}\label{eq: E3 2-dim d1}\begin{array}{l}
\Delta_3=-m^{12}\left(2m^3\mu_1
a+m^2\mu_1^2+m^4a^2+2m\mu_1+1+2m^2a-4m^3\right).\end{array}\end{eqnarray*}
the singular fibers are $I_6^*+I_2^*+4I_1$ and $MW=\Z/2\Z$; \item
$\mathcal{E}_4$ induced by the conic bundle $x^2+axy+dxz+yz$,
whose discriminant is
\begin{eqnarray*}\label{eq: E4 2-dim d1}\begin{array}{l}
\Delta_4:=-d^2\left(-1+da\right)^2\left(d^4-2d^3\mu_1+\mu_1^2d^2+2d^2a-4d-2a\mu_1
d+4\mu_1+a^2\right)\end{array}\end{eqnarray*} the singular fibers
are $I_{10}^*+2I_2+4I_1$ and $MW=\Z/2\Z$;

\item $\mathcal{E}_5$ induced by the conic bundle $bx^2+axy+yz$,
whose discriminant is
\begin{eqnarray*}\label{eq: E5 2-dim d1}\begin{array}{l}
\Delta_5:=b^9(-1+b)^2\left(4ba^3-4a^3-a^2b\mu_1^2-18b\mu_1
a+18b^2\mu_1 a-4b^2\mu_1^3-54b^2+27b^3+27b\right)
\end{array}\end{eqnarray*}
the singular fibers are $II^*+III^*+I_2+3I_1$ and $MW=\{1\}$;
\end{enumerate}
\textbf{}
\\
{\bf Specializations of the rational elliptic surface.} If $a=0$,
the rational elliptic surface $R_a$ becomes the rigid rational
elliptic surface with reducible fibers $I_9+3I_1$, so we go back
to our previous case.

From now on we assume that $a\neq 0$. If
$$a\in\left\{-\frac{3}{8}\sqrt[3]{2},
\frac{3\left(\sqrt[3]{2}-i\sqrt{3}\sqrt[3]{2}\right)}{16},\frac{3\left(\sqrt[3]{2}+i\sqrt{3}\sqrt[3]{2}\right)}{16}\right\},$$
then the discriminant $\Delta$ has a multiple zero and in
particular if $a=-\frac{3}{8}\sqrt[3]{2}$ (resp.
$a=\frac{3\left(\sqrt[3]{2}\pm i\sqrt{3}\sqrt[3]{2}\right)}{16}$),
then the fibration has a fiber of type $II$ in
$\mu_1=-\frac{3}{2}\sqrt[3]{4}$ (resp. $\mu_1=\frac{\mp
3i\sqrt[3]{4}\sqrt{3}+\sqrt[3]{4}}{4}$). The other fibers are
unchanged and they are $I_8+2I_1$.
\\
\\
{\bf Specialization of the surface $X_{a,\mu_1}$}
Let us assume that $a\neq 0$ and
$a\not\in\left\{-\frac{3}{8}\sqrt[3]{2},
\frac{3\left(\sqrt[3]{2}\pm
i\sqrt{3}\sqrt[3]{2}\right)}{16}\right\}$. If $\mu_1$ is chosen to
be one of the solutions of
$$
\Delta=-\frac{1}{16}\mu^4a+\frac{1}{16}\mu^3+\frac{1}{2}\mu^2a^2-\frac{9}{4}\mu
a-a^3+\frac{27}{16}=0,$$ then the K3 surface $X_{a,\mu_1}$ is
constructed as a double cover of $R_a$, branched over a fiber of
type $I_8$ and a fiber of type $I_1$. So $\mathcal{E}_1$ is an
elliptic fibration with fibers $I_{16}+I_2+6I_1$ and $MW=\Z$ and
the K3 surface $X_{a,\mu_1}$ admits an involution acting trivially
on the N\'eron--Severi group whose fixed locus consists of 9
rational curves. For this specialization $T_{X_{\mu_1}}\simeq
U(2)\oplus \langle 2\rangle$ (because of the presence of such an
automorphism). Since the degree of $\Delta$ is 4, this happens for
4 values of $\mu_1$.

If $\mu_1=1/a$ one observes that the discriminant of the elliptic
fibration $\mathcal{E}_2$ becomes $m^{10}(a+m)^3(4m^2a^2-a-m)$.
The solution $m=-a$ is now a solution of order 3 (and not only of
order 2 as it is for general values of $\mu_1$). Moreover now
$m=-a$ is also a solution of $A_2$ and $B_2$ with multiplicities 1
and 2 respectively. This implies that the elliptic fibration
$\mathcal{E}_2$ has now a fiber of type $III$ and the singular
fibers of this fibration are $III^*+I_4^*+III+3I_1$. This has no
effect on the rank of the N\'eron Severi group.

One can consider the discriminant of the other elliptic
fibrations, obtaining 5 values of $\mu_1$ for which the Picard
number of $X_{\mu_1}$ increases, 4 are solutions of $\Delta=0$ and
one is $\frac{-a^2}{4}$. If $\mu_1=\frac{-a^2}{4}$, then the
discriminant of $\mathcal{E}_5$ has a zero with multiplicity three
(which is of multiplicity 2 if $\mu_1\neq \frac{-a^2}{4}$), so the
trivial lattice of $\mathcal{E}_5$ becomes $U\oplus E_8\oplus
E_7\oplus A_2$ and since there can not be torsion section we
conclude that this is also the N\'eron-Severi group. So the
transcendental lattice is $T_X\simeq A_2(-1)\oplus \langle
-2\rangle$.

We put in a Table the information about the specializations due to
$\mu_1$ which changes the N\'eron--Severi of the surface:
$$
\begin{array}{|c|c|c|}
\hline
\mu_1&T_{X_{\mu_1}}&\mbox{ fibrations }\\
\hline
\mbox{ solutions }\Delta=0&\langle2\rangle^2\oplus\langle-2\rangle &\begin{array}{|c|c|c|}\mathcal{E}_i&\mbox{ fibers }&MW\\
\hline
\mathcal{E}_1&I_{16}+I_2+6I_1&\Z\\
\mathcal{E}_2&III^*+I_4^*+2I_2+I_1
&\Z/2\Z\\
\mathcal{E}_3&I_6^*+I_2^*+I_2+2I_1&\Z/2\Z\\
\mathcal{E}_4&I_{10}^*+3I_2+2I_1&\Z/2\Z\\
\mathcal{E}_5&II^*+III^*+2I_2+I_1&\{1\}\\
\end{array}\\
\hline \frac{-a^2}{4}&A_2(-1)\oplus\langle 2\rangle&
\begin{array}{|c|c|c|}\mathcal{E}_i&\mbox{ fibers }&MW\\
\hline
\mathcal{E}_1&I_{16}+8I_1&\Z\times \Z\\
\mathcal{E}_2&III^*+I_4^*+I_2+3I_1&\Z\times \Z/2\Z\\
\mathcal{E}_3&I_6^*+I_2^*+4I_1&\Z\times \Z/2\Z\\
\mathcal{E}_4&I_{10}^*+I_3+I_2+3I_1&\Z/2\Z\\
\mathcal{E}_5& II^*+III^*+I_3+2I_1&\{1\}
\end{array}\\
\hline
\end{array}
$$
From the previous table it is clear that for $a=0$ one obtains the
previous specializations to K3 surfaces with transcendental
lattice either $\langle 2\rangle^2$ or $A_2(-1)$.

Let us now assume that $a=-\frac{3}{8}\sqrt[3]{2}$ (the cases $a=
\frac{3\left(\sqrt[3]{2}\pm i\sqrt{3}\sqrt[3]{2}\right)}{16}$ are
analogous). In this case the rational elliptic fibration has
fibers of types $I_8+II+2I_1$. If the fiber over $\mu_1$ is smooth
nothing changes in the N\'eron--Severi of the K3 surface with
respect to the general case. If the fiber over $\mu_1$ is $I_1$
nothing changes in the N\'eron--Severi of the K3 surface with
respect to the choice of a different value of $a$, and $\mu_1$ to
be a solution of $\Delta=0$. But if the fiber over $\mu_1$ is of
type $II$, then something changes. In particular we have the
following Table

$$
\begin{array}{|c|c|c|}
\hline
&a=-\frac{3}{8}\sqrt[3]{2}&\\
\hline
\mu_1&T_{X_{\mu_1}}&\mbox{ fibrations }\\
\hline
-\frac{3}{2}\sqrt[3]{4}&\langle2\rangle\oplus\langle 6\rangle&\begin{array}{|c|c|c|}\mathcal{E}_i&\mbox{ fibers }&MW\\
\hline
\mathcal{E}_1&I_{16}+IV+4I_1&\Z\\
\mathcal{E}_2&III^*+I_4^*+I_3+I_2
&\Z/2\Z\\
\mathcal{E}_3&I_6^*+I_2^*+I_3+I_1&\Z/2\Z\\
\mathcal{E}_4&I_{10}^*+I_3+2I_2+I_1&\Z/2\Z\\
\mathcal{E}_5&II^*+III^*+I_2+I_3&\{1\}\\
\end{array}\\
\hline
\end{array}
$$
\begin{rem}{\rm In \cite{wine1}, the authors classify all the elliptic fibrations on the unique K3 surface with transcendental lattice $\langle 2\rangle\oplus\langle 6\rangle$. Here we recover some of the elliptic fibrations on the mentioned K3 surface, in particular the ones numbered as 17,12,10, 5 and 2 in the aforementioned article. Their equations are easily obtained here, just specializing the ones of $\mathcal{E}_i$ on $X$. One observes that in the case of the fibration 17, which corresponds here to $\mathcal{E}_1$, the lattice $A_2$ in the trivial lattice corresponds not to a fiber of type $I_3$ (as in general happens) but to a fiber of type $IV$.}\end{rem}

\subsubsection{$\delta=0$}
We consider the 1-dimensional family of rational elliptic surfaces, parametrized by $a$, associated to the pencil of cubics:
$$x^2y+z^2x+xy^2+az^2y+\mu xyz$$
The Weierstrass equation of the rational elliptic fibration
$\mathcal{E}_R:R_a\ra\mathbb{P}^1_{\mu}$ is
\begin{eqnarray*}\begin{array}{ll}Y^2=&X^3+\left(-\frac{1}{48}\mu^4+\frac{1}{6}\mu^2a+\frac{1}{6}\mu^2-\frac{1}{3}a^2+\frac{1}{3}a-1/3\right)X\\&+\frac{1}{864}\left(4a+4-\mu^2\right)\left(16a^2-8\mu^2a-40a+16+\mu^4-8\mu^2\right)\end{array}\end{eqnarray*}
and the discriminant is
$$\Delta:=-\frac{1}{16}a^2\left(4a-4-\mu^2+4\mu\right)\left(4a-4-\mu^2-4\mu\right).$$

To obtain the K3 surface $X_{a,\mu_1}$ with the elliptic fibration
$\mathcal{E}_X:X_{a,\mu_1}\ra\mathbb{P}^1_{\tau}$ we apply the
base change $\mu=\tau^2+\mu_1$ which is branched at infinity and
$\mu=\mu_1$.

On $X_{a,\mu_1}$ there are 6 elliptic fibrations.

They are:\begin{enumerate}\item $\mathcal{E}_1$ induced by the
elliptic fibration $\mathcal{E}_R$ on $R_a$ whose discriminant is
\begin{eqnarray*}\label{eq: E1 2-dim d2}
\begin{array}{l}
\Delta_1:=-\frac{1}{16}a^2\left(4a-4-(\tau^2+\mu_1)^2+4(\tau^2+\mu_1)\right)\left(4a-4-(\tau^2+\mu_1)^2-4(\tau^2+\mu_1)\right)
\end{array}
\end{eqnarray*}
the singular fibers are $I_{16}+8I_1$ and $MW=\Z\times \Z/2\Z$; \item
$\mathcal{E}_2$ induced by the conic bundle $x-my$, whose
discriminant is
\begin{eqnarray*}\label{eq: E2 2-dim d2} \begin{array}{l}
\Delta_2:=m^9(1+m)^2(a+m)^2(4ma+4a+4m^2+4m-m\mu_1^2)\end{array}\end{eqnarray*}
the singular fibers are $2III^*+2I_2+2I_1$ and $MW=\Z/2\Z$; \item
$\mathcal{E}_3$ induced by the conic bundle $y-mz$, whose
discriminant is
\begin{eqnarray*}\label{eq: E3 2-dim d2}\begin{array}{l}
\Delta_3:=a^2m^{10}(-m^4-2m^3\mu_1-m^2\mu_1^2-2m^2+4am^2-2m\mu_1-1)
\end{array}\end{eqnarray*}
the singular fibers are $2I_4^*+4I_1$ and $MW=\Z/2\Z$; \item
$\mathcal{E}_4$ induced by the conic bundle $x+y+mz$, whose
discriminant is
\begin{eqnarray*}\label{eq: E4 2-dim d2} \begin{array}{l}
\Delta_4:=-a^2m^6\left(m^4-2m^3\mu_1+m^2\mu_1^2+2am^2+2m^2-2m\mu_1-2m\mu_1a+1-2a+a^2\right)
\end{array}\end{eqnarray*}
the singular fibers are $I_8^*+I_0^*+4I_1$ and $MW=\Z/2\Z$; \item
$\mathcal{E}_5$ induced by the conic bundle $xy+y^2+dxz $, whose
discriminant is
\begin{eqnarray*}\label{eq: E5 2-dim d2} \begin{array}{l}
\Delta_5:=d^{10}\left(4d^2-8da\mu_1-d^2\mu_1^2+2a+2d\mu_1-6d^2a^2+4a^3d^2+18d^3a\mu_1+\right.\\
\left.
10d^2a\mu_1^2-d^2a^2\mu_1^2+2a^2d\mu_1-4d^3a\mu_1^3-6d^2a-1+18d^3a^2\mu_1+27d^4a^2-a^2\right)
\end{array}\end{eqnarray*}
the singular fibers are $II^*+I_4^*+4I_1$ and $MW=\{1\}$;

\item $\mathcal{E}_6$ induced by the conic bundle
$x^2y+2xy^2+y^3+z^2x+hyz(x+y)$, whose discriminant is
\begin{eqnarray*}\label{eq: E5 2-dim d2} \begin{array}{l}
\Delta_6:=-4096(4a+a^2h^6+2a^3h^4+a^4h^2+8h^5\mu_1a+12a^2-2\mu_1h^3+12a^3+4a^4-27\mu_1^2a^2-2h^6a+\\+h^4\mu_1^2-2h^5\mu_1+h^6+2h^4+6a^2h^3\mu_1-6a^3h\mu_1
+6\mu_1h^3a-9a^2h^2+42a^2h\mu_1+10a^3h^2+10h^2a+\\
h^2+a^2h^4\mu_1^2-2a^2h^5\mu_1-2a^3h^3\mu_1-6ah^2\mu_1^2-6a^2h^2\mu_1^2+4ah^3\mu_1^3-10ah^4\mu_1^2-6ah\mu_1)a^2\\
\end{array}\end{eqnarray*}
the singular fibers are $I_{12}^*+6I_1$ and $MW=\{1\}$. In order
to compute the Weierstrass equation as in the algorithm we perform
a change of coordinates sending the point $(1:-1:0)$ to $(0:1:0)$.
\end{enumerate}

We observe that $a\neq 0$, otherwise the pencil of cubics defining
the rational elliptic surface does not contain smooth fibers.

If $\mu_1=2(\pm 1\pm\sqrt{a})$, then the branch fibers of the
cover $X\ra R$ are $I_8$ and $I_1$.

If $\mu_1=0$ the elliptic fibration $\mathcal{E}_1$ becomes
$Y^2=X^3+\frac{1}{3}m^3(1+m)(a+m)X$ which is an isotrivial
fibration whose general fiber is isometric to the elliptic curve
with complex multiplications of order 4. In this case there is an
extra automorphism of order 4, but the N\'eron--Severi does not
change. The order 4 automorphism is purely non--symplectic and the
presence of such an automorphism reduces the dimension of the
family of K3 surface from 2 to 1, see \cite{AS}. The
specializations of the elliptic fibrations $\mathcal{E}_i$ on
$X_{a,\mu_1}$ for specific values of $a$ and $\mu_1$ are
summarized in the following table.

$$
\begin{array}{|c|c|c|}
\hline
\mbox{ specialization }&T_{X_{\mu_1}}&\mbox{ fibrations }\\
\hline
\mu_1=2(\pm 1\pm\sqrt{a})&\langle2\rangle^2\oplus\langle-2\rangle &\begin{array}{|c|c|c|}\mathcal{E}_i&\mbox{ fibers }&MW\\
\hline
\mathcal{E}_1&I_{16}+I_2+6I_1&\Z\times \Z/2\Z\\
\mathcal{E}_2&2III^*+3I_2
&\Z/2\Z\\
\mathcal{E}_3&2I_4^*+I_2+2I_1&\Z/2\Z\\
\mathcal{E}_4&I_8^*+I_0^*+I_2+2I_1&\Z/2\Z\\
\mathcal{E}_5&II^*+I_4^*+I_2+2I_1&\{1\}\\
\mathcal{E}_6&I_{12}^*+I_2+4I_1&\{1\}
\end{array}\\
\hline a=1&U\oplus\langle 4\rangle&
\begin{array}{|c|c|c|}\mathcal{E}_i&\mbox{ fibers }&MW\\
\hline
\mathcal{E}_1&I_{16}+8I_1&\Z\times \Z\times \Z/2\Z\\
\mathcal{E}_2&2III^*+I_4+2I_1
&\Z/2\Z\\
\mathcal{E}_3&2I_4^*+4I_1&\Z/2\Z\times \Z\\
\mathcal{E}_4&I_8^*+I_1^*+3I_1&\Z/2\Z\\
\mathcal{E}_5&II^*+I_5^*+3I_1&\{1\}\\
\mathcal{E}_6&I_{12}^*+6I_1&\Z
\end{array}\\
\hline a=1,\ \mu_1=0&\langle 2\rangle^2&
\begin{array}{|c|c|c|}\mathcal{E}_i&\mbox{ fibers }&MW\\
\hline
\mathcal{E}_1&I_ {16}+I_4+4I_1&\Z/4\Z\\
\mathcal{E}_2&2III^*+I_0^*
&\Z/2\Z\\
\mathcal{E}_3&2I_4^*+2I_2&(\Z/2\Z)^2\\
\mathcal{E}_4&I_8^*+I_2^*+2I_1&\Z/2\Z\\
\mathcal{E}_5&II^*+I_6^*+2I_1&\{1\}\\
\mathcal{E}_6&I_{14}^*+4I_1&\Z
\end{array}\\
\hline a=1,\ \mu_1=\pm 4&\langle 2\rangle\oplus\langle 4\rangle&
\begin{array}{|c|c|c|}\mathcal{E}_i&\mbox{ fibers }&MW\\
\hline
\mathcal{E}_1&I_{16}+3I_2+2I_1&\Z/4\Z\\
\mathcal{E}_2&2III^*+I_2+I_4&\Z/2\Z\\
\mathcal{E}_3&2I_4^*+I_2+2I_1&\Z\times \Z/2\Z\\
\mathcal{E}_4&I_8^*+I_1^*+I_2+I_1&\Z/2\Z\\
\mathcal{E}_5&II^*+I_5^*+I_2+I_1&\{1\}\\
\mathcal{E}_6&I_{13}^*+I_2+3I_1&\{1\}
\end{array}\\
\hline
\end{array}
$$
\begin{rem}{\rm The elliptic fibrations on the K3 surface with transcendental lattice $\langle 2\rangle\oplus \langle 4\rangle$ are classified in \cite{BL}, where their Weierstrass equations are also given. The fibrations $\mathcal{E}_i$ $i=1,\ldots, 6$ corresponds to the fibrations described lines 18, 12, 14, 22, 25, 23 of the last table in \cite{BL} respectively. In our context the equations are just obtained by the application of the algorithm with the assumption $a=1$, $\mu_1=\pm 4$.}\end{rem}

\section{Higher genus curve in the fixed locus of $\iota$}\label{sec: higher genus}
In what follows, we assume that $X$ is a K3 surface and $\iota$ is a
non-symplectic involution on $X$ whose fixed locus contains a
curve $C$ of genus higher than 1 and $k$ rational curves. We assume that $\iota^*$ acts as the identity on the N\'eron--Severi
group.

\begin{proposition}\label{prop: higher genus}
Let $(X,\iota)$ be as in Assumption \ref{cond: X iota} and
$\mathcal{E}:X\ra\mathbb{P}^1$ an elliptic fibration on $X$. Let
$C$ be a genus $g>1$ curve fixed by $\iota$.

Then the elliptic fibrations on $X$, necessarly of type 1 with
respect to $\iota$, are listed given in the Appendix.

Moreover each elliptic fibration induces a rational fibration on $X/\iota$ (and thus a pencil of rational curves on any birational model of $X/\iota$).
\end{proposition}

\proof By Theorem \ref{prop: when iota is trivial on the basis},
each elliptic fibration on $X$ is of type 1, so by Proposition
\ref{prop: possible fibers} one has the list of the reducible
fibers which can appear on elliptic fibrations on $X$. Moreover,
by Proposition \ref{prop: mordell weil}, the Mordell-Weil group of
any fibration on $X$ is contained in $\Z/2\Z$ and it is trivial if
$g>3$. This allows us to produce the list of elliptic fibrations
exactly as in Proposition \ref{prop: possible fibrations}. Since
the involution $\iota$ is the cover involution, it induces a
fibration on $X/\iota$ whose fiber are the quotient of the fiber
of the fibration on $X$ by the restriction of $\iota$. Since
$\iota$ is the elliptic involution on the fiber, the fibration on
$X/\iota$ has rational fibers.

All the fibrations listed has Mordell Weil rank which is 0, then
they exist since they appear in the list in \cite{Shim}.
\endproof

\begin{rem}{\rm The classifications of the elliptic fibrations given in Section \ref{sec: appendix} can be also obtained by Nishiyama's method, considering the following lattices:
$$\begin{array}{|c|c||c|c|}
\hline g,k&T&g,k,\delta&T\\
\hline
g=2,\ k=1,\ldots,9&E_8\oplus A_1^{9-k}&g=2,\ k=5,\ \delta=0&D_8\oplus D_4\\
\hline g=3,\ k=1,\ldots,6& E_8\oplus D_4\oplus A_1^{6-k}& g=3,\
k=2,\ \delta=0&D_8\oplus D_4^2\\
\hline g=4,\ k=1,\ldots,5&E_8\oplus D_6\oplus A_1^{5-k}&g=4,\
k=3,\
\delta=0& E_8\oplus D_4\oplus D_4\\
\hline
g=5,\ k=1,\ldots 5& E_8\oplus E_7\oplus A_1^{5-k}& g=5,\
k=4,\
\delta=0&D_8\oplus D_8\\
\hline
g=6,\ k=1,\ldots 5&E_8\oplus E_8\oplus A_1^{5-k}&g=7,\
k=1,2&E_8\oplus E_8\oplus D_4 \oplus A_1^{2-k}\\
\hline g=8,\ k=1&E_8\oplus E_8\oplus D_6&g=9,\ k=1&E_8\oplus
E_8\oplus
E_7\\
\hline g=10,\ k=1& E_8\oplus E_8\oplus E_8\\
\hline
\end{array}$$
}\end{rem}

The K3 surfaces $X$ as in Proposition \ref{prop: higher genus}
admit two different very natural geometric descriptions: one is as
double cover of a minimal model of $X/\iota$ (which is rational
and in some cases $\mathbb{P}^2$), the other is
$\varphi_{|C|}(X)\subset \mathbb{P}^{g}$, if $g>2$.

\subsection{The K3 surfaces $X$ as double covers of $\mathbb{P}^2$}

Let us consider the double cover of $\mathbb{P}^2$ branched on a
possibly reducible sextic. The minimal model of this double cover
is a K3 surface $X$, admitting a non-symplectic involution $\iota$
(the cover involution) and for sufficiently generic choices of the
sextic, $\iota$ acts trivially on the N\'eron--Severi group of
$X$. Hence the N\'eron--Severi group of $X$ can be deduced by the
Nikulin classification of the non-symplectic involutions and
depends on the number and the genus of the components of the
branch sextic (whose normalization is isomorphic to the fixed
locus of $\iota$). On the other hand the choice of the
N\'eron--Severi group of a K3 surface which admits a
non-symplectic involution acting trivially on the N\'eron--Severi
group, determines a K3 surface which satisfies the hypothesis of
Proposition \ref{prop: higher genus}, if the fixed curve with
highest genus has genus at least 2.

By Proposition \ref{prop: higher genus}, one obtains that the elliptic fibrations on $X$ are induced by ``generalized conic bundles", i.e. by pencils of rational curves passing through a certain numbers of singular points of the branch sextic. This allows to reproduce all the computations did in case $g=1$ also in these highest genus cases, if one is able to describe explicitly the branch sextic.

If $g=6$ and $k=1$, then the K3 surface $X$ is generic in the
unique family of K3 surfaces admitting an involution $\iota$
fixing one rational curve and a curve of genus 6. It can be
realized as minimal model of the double cover of $\mathbb{P}^2$
branched on a line $l$ and an smooth quintic $q$.

Similarly, different values of $g$ and $k$ can be obtained modifying the branch curve. So some of the surfaces in Proposition \ref{prop: higher genus} can be realized as minimal model of a double cover of $\mathbb{P}^2$ branched on a sextic $s$ as in the table:
$$\begin{array}{|l|c|c|c|}
\hline
\mbox{ sextic $s$}&g&k\\
\hline
\mbox{line $l$+quintic $q$,\ $q$ has $0\leq \alpha\leq 4$ nodes}&6-\alpha&1\\
\hline
\mbox{line $l$+quintic $q$,\ $q$ has $\alpha=1,2$ nodes, $m$ through the node of $q$}&6-\alpha&2\\
\hline
\mbox{line $l$+line $m$+quartic $q$,\ $q$ has $\alpha=0,1$ node}&3-\alpha&2\\
\hline
\mbox{line $l$+line $m$+quartic $q$,\ $q$ has $\alpha=0,1$ node, $l\cap m\cap q\neq \emptyset$}&3-\alpha&3\\
\hline
\begin{array}{l}\mbox{line $l$+line $m$+quartic $q$,\ $q$ has $1$ node, $l\cap m\cap q\neq \emptyset$,}\\\mbox{$m$ through the node of $q$}\end{array}&2&4\\
\hline
\end{array}
$$
We observe that if the curve $q$ has a node, then the pencil of lines through this node induces an elliptic fibration on the K3 surface $X$, double cover of $\mathbb{P}^2$ branched on the sextic $s$. This provides a geometric description, for example, of the elliptic fibrations obtained listed in the Appendix for $g=2,\ldots 6$, $k=1$.

\subsection{The model associated to $|C|$, $g=2,3$}
The curve $C$ fixed by $\iota$ is a smooth irreducible curve on $X$, so its linear system defines a map to a projective space, in particular to $\mathbb{P}^g$. For low values of $g$, this model is quite clear and sometime also provide an equation for $X$. Here we collect some results on this.
\begin{proposition}\label{prop: high genus=2}
Let $(X,\iota)$ be as in Assumption \ref{cond: X iota},
$\mathcal{E}:X\ra\mathbb{P}^1$ an elliptic fibration on $X$ and
$C\in Fix_{\iota}$ a curve of genus 2. Then, up to a choice of the
coordinates of $\mathbb{P}^2$, $\varphi_{|C|}:X\ra \mathbb{P}^2$
is a double cover of $\mathbb{P}^2$ branched over the sextic
\begin{equation*}
x_0^4f_2(x_1:x_2)+x_0^2f_4(x_1:x_2)+f_6(x_1:x_2), \end{equation*}
where $f_i\in\C_i[x_1:x_2]$,  and $\iota$ is induced on $X$ by
$\iota_{|\mathbb{P}^2}:(x_0,x_1,x_2)\ra(-x_0,x_1,x_2)$. Denote by
$\alpha$ the cover involution. The following hold:
\begin{enumerate}
\item[i)] If $C$ is a trisection of $\mathcal{E}$, the fibers of
the elliptic fibration $\mathcal{E}$ are mapped to curves of
degree 3 in $\mathbb{P}^2$ which passes through $(1:0:0)$ and
intersects the branch sextic in every intersection point with even
multiplicity, and $\alpha$ is an involution of the basis of
$\mathcal{E}$. \item[ii)] If $C$ is a bisection of $\mathcal{E}$,
the fibers of the elliptic fibration $\mathcal{E}$ are mapped to
lines in $\mathbb{P}^2$ which passes through $(1:0:0)$. Let
$\sigma$ be the symplectic involution which is the translation by
the 2-torsion section of $\mathcal{E}$, then $\alpha=\iota\circ
\sigma$.
\end{enumerate}
\end{proposition}

\proof Since the curve $C$ is a genus 2 curve, the map
$\varphi_{|C|}:X\ra\mathbb{P}^2$ is a double cover. The cover
involution $\alpha$ is not $\iota$, since $\alpha$ fixes the
branch locus of the cover $X\ra\mathbb{P}^2$ and $\iota$ fixes the
pullback of a generic the hyperplane section, i.e. $C$. The
automorphism $\iota$ preserves the map $X\ra\mathbb{P}^2$, and
fixes the pull-back of a line, so it descends to an automorphism
$\iota_{\mathbb{P}^2}$ of $\mathbb{P}^2$ which fixes the branch
locus of the cover $X\ra\mathbb{P}^2$. Up to a choice of
coordinates we can assume that $\iota_{\mathbb{P}^2}$ is
$(x_0:x_1:x_2)\mapsto (-x_0:x_1:x_2)$. This fixes the line $x_0=0$
and  the point $(1:0:0)$. The sextics invariant for $\iota_{\mathbb{P}^2}$ have equation
$ax_0^6+x_0^4f_2(x_1:x_2)+x_0^2f_4(x_1:x_2)+f_6(x_1:x_2)$.

Since $X$ admits an elliptic fibration, $\iota$ fixes not only $C$
but at least on rational section, i.e. a rational curve $R_1$
orthogonal to $C$. So $\varphi_{|C|}$ contracts at least one
rational curve, fixed by $\iota$. The branch sextic of
$\varphi_{|C|}:X\ra\mathbb{P}^2$ has a singularity in the fixed
point of $\iota_{\mathbb{P}^2}$. Hence $a=0$ and the branch sextic
is $x_0^4f_2(x_1:x_2)+x_0^2f_4(x_1:x_2)+f_6(x_1:x_2)$, the inverse
image of $x_0=0$ on $X$ is exactly the smooth genus 2 curve $C$,
the curve(s) resolving the singularity of the branch sextic in
$(1:0:0)$ are rational curve(s) possibly fixed by $\iota$.

If $C$ is a trisection, denoted by $F$ the class of the fiber of
$\mathcal{E}$, $C \cdot F=3$, hence the image of each curve in $F$ has
degree 3 in $\mathbb{P}^2$ and the curves in $|\varphi_{|C|}(F)|$
splits in the double cover (otherwise $C\cdot F$ should be an even
number). Hence $\alpha$ is the involution which switches the two
disjoint curves in the inverse image of the curves in
$|\varphi_{|C|}(F)|$. So it is an involution which preserves the
class of the fiber but switches pairs of fibers of $\mathcal{E}$,
and it is an involution of the basis of $\mathcal{E}$. Moreover,
$F\cdot R_1=1$, so the curves in $|\varphi_{|C|}(F)|$ pass through
$\varphi_{|C|}(R_1)=(1:0:0)$.

If $C$ is a bisection, $C \cdot F=2$, and there are two possibilities:
either $\varphi_{|C|}(F)$ is a curve of degree 2 in $\mathbb{P}^2$
and it splits in the double cover, or  $\varphi_{|C|}(F)$ is a
curve of degree 1 in $\mathbb{P}^2$ which does not split. If a
conic splits in the double cover, its inverse image consists of
two copies of a rational curve. This can not be the case, since $|F|$ is 1-dimensional system of genus 1 curves. So $\varphi_{|C|}(F)$ is
a line (which does not split in the double cover).
The lines in $|\varphi_{|C|}(F)|$ pass through the node, because $FR_1=1$.

Since $C$ is a bisection of $\mathcal{E}$ passing through some of
the 2-torsion points, there exists also a 2-torsion section $T$ in
$\mathcal{E}$, which is fixed by $\iota$ and contracted by $\varphi_{|C|}$ since
$CT=0$. In particular the singularity $(1:0:0)$ in the branch
sextic is worst than a simple node. Due to the symmmetry
$(x_0:x_1:x_2)\ra(-x_0:x_1:x_2)$, the inverse image on $X$ of
point $(1:0:0)$ is not simply an $A_2$ configuration of rational
curves but it is at least an $A_3$ configuration of rational
curves. The two rational curves at the extreme of the trees of
rational curves of these $A_3$ are switched by $\alpha$ and they
are sections of the elliptic fibration $\mathcal{E}$. These two curves correspond to the zero section and
the 2-torsion section of $\mathcal{E}$, both are fixed by $\iota$. The
involutions $\alpha$ and $\iota$ preserve each fiber of the
fibration, so the same do their composition $\alpha\circ
\iota=\sigma$, which is then a symplectic involution on
$\mathcal{E}$, preserving each fiber and mapping the zero section
to the 2-torsion section. This implies that $\sigma$ is a
translation by the 2-torsion section.
\endproof
\begin{rem}{\rm
By Proposition \ref{prop: high genus=2}, it follows that the K3
surface $X$ as in the hypothesis admits not only a non-sympelctic
involution $\iota$, but also a symplectic involution $\sigma$.
These K3 surfaces were studied in \cite{GS sympl and non sympl},
where numeric conditions on $g$ and $k$ which imply that K3
surfaces with a non-symplectic involution necessarily admit also a
symplectic one are established. The conditions in
\cite[Proposition 3.1]{GS sympl and non sympl} imply (with just
one exception) that the fixed curves have genus at most 2. Hence,
Proposition \ref{prop: high genus=2} gives a geometric
interpretation of the involutions considered in \cite{GS sympl and
non sympl}}\end{rem}

Let us now briefly discuss the case $g=3$, naturally associated to a map to $\mathbb{P}^3$. We observe that there is just one case in which this map is $2:1$ onto the image, and not $1:1$. This is exactly the case mentioned in the previous remark, where one has a curve of genus 3 in the fixed locus of an involution $\iota$ acting trivially on the N\'eron--Severi group, but it is still true that every K3 surfaces with this property also admit a symplectic involution.

\begin{proposition}
Let $(X,\iota)$ be as in Assumption \ref{cond: X iota},
$\mathcal{E}:X\ra\mathbb{P}^1$ an elliptic fibration on $X$ and
$C\in Fix_{\iota}$ a curve of genus 3.

If $C$ is a trisection, then $\varphi_{|C|}:X\ra\mathbb{P}^3$ is a
generically $1:1$ map to a singular quartic with equation
$$x_0^2f_2(x_1:x_2:x_3)+f_4(x_1:x_2:x_4).$$
The eight lines connecting the singular point $(1:0:0:0)$ to the eight points $f_2\cap f_4\subset\mathbb{P}^2_{(x_1:x_2:x_3)}$ are contained in $\varphi_{|C|}(X)$.

If $C$ is a bisection, then $\varphi_{|C|}:X\ra\mathbb{P}^3$ is a
generically $2:1$ map to a quadric.

\end{proposition}
\proof Let us now assume that $X$ admits an elliptic fibration
$\mathcal{E}$ and that $C$ is a trisection of this fibration. Then
$C$ in not hyperelliptic and by \cite{SD} the map
$\varphi_{|C|}:X\ra\mathbb{P}^3$ exhibits $\varphi_{|C|}(X)$ as
quartic hypersurface in $\mathbb{P}^3$. Denote by $F$ and $O$
respectively the fiber and the zero section of $\mathcal{E}$. We have
$C\cdot O=0$, so the curve $O$ is contracted by $\varphi_{|C|}$ and thus the
quartic has a node. Moreover, the involution $\iota$ descends to
an involution of $\mathbb{P}^3$, fixing the node. Up to a choice
of coordinates, one can assume that the node is $(1:0:0:0)$ and the
involution induced by $\iota$ on $\mathbb{P}^3$ is
$(x_0:x_1:x_2:x_3)\mapsto (-x_0:x_1:x_2:x_3)$. Since $F \cdot C=3$, the
fibers of $\mathcal{E}$ are mapped to curves of degree 3 in
$\mathcal{E}$, passing through the node (since $F \cdot O=1$).
By the equation of the quartic one immediately checks that there are eight lines contained in $\varphi_{|C|}(X)$ passing through $O$. Let $l$ be one of these lines and $L$ be the class in $NS(X)$ corresponding to the strict transform of $l$ after blowing up the node. We have $L^2=-2$, $L \cdot C=1$, $L \cdot O=1$. The pencil of hyperplanes through $l$ cuts on $\varphi_{|C|}(X)$ a pencil of genus 1 curves, passing through the node, so it induces on $X$ an elliptic fibration, whose class is $C-L$. Generically this pencil has 7 reducible fibers, corresponding to the hyperplane through $l$ which contains another line of $\varphi_{|C|}(X)$. Indeed generically the reducible fibers of $\mathcal{E}$ are  7 fibers of type
$I_2$ (see Table \ref{table of admissible fibers  g=3}, $k=1$).

If $C$ is a bisection, it is a hyperelliptic curve and by \cite{SD} the map $\varphi_{|C|}:X\ra\mathbb{P}^3$ is $2:1$ to a quadric.
\endproof

\section{Appendix}\label{sec: appendix}
In this section we list the elliptic fibrations on a K3 surface
admitting an involution $\iota$ which acts trivially on the
N\'eron--Severi group and such that the highest genus $g$ of a
fixed curve is greater than 1. The lists are obtained according to
the results of Theorem \ref{prop: when iota is trivial on the
basis} and Propositions \ref{prop: mordell weil} and \ref{prop:
possible fibers} similarly to what is done in the proof of
Proposition \ref{prop: possible fibrations}. Also the notation is
the same as in Propostion \ref{prop: possible fibrations}. The
fibrations on K3 surfaces admitting an involution as $\iota$ such
that $g=2$ are given in Table \ref{table of admissible fibers
g=2}; the ones such that $g=3$ in Table \ref{table of admissible
fibers g=3};
 the ones such that $g=4$ in Table \ref{table of admissible fibers g=4}
; the ones such that $g=5$ in Table \ref{table of admissible fibers g=5}
; the ones such that $g=6$ in Table \ref{table of admissible fibers g=6}
; the ones such that $g=7$ in Table \ref{table of admissible fibers g=7}
; the ones such that $g=8$ in Table \ref{table of admissible fibers g=8+b}, $b=0$
; the ones such that $g=9$ in Table \ref{table of admissible fibers g=8+b}, $b=1$
; the ones such that $g=10$ in Table \ref{table of admissible fibers g=8+b}, $b=2$.

\begin{longtable}{l}\caption{}\label{table of admissible fibers g=2}\\ \hline $g=2, k=9, r=18, a=0$, ($\delta$=0)\\
$
\begin{array}{|c|c|c|c|c|}
\hline
\mbox{trivial lattice}&16=\sum c_i+\rk(MW)&9=k=\sum s_i+\#\mbox{sections}&MW(\mathcal{E})\\
\hline
U\oplus E_8\oplus E_8&8+8&4+4+1&\{1\}\\
\hline
U\oplus D_{16}&16& 7+2&\Z/2\Z\\
\hline
\end{array}$\\
 $ $\\
\hline
\mbox{$g=2,\ k=8, r=17,\ a=1$,\ ($\delta$=1)}\\
$\begin{array}{|c|c|c|c|}
\hline
\mbox{trivial lattice}&15=\sum c_i+\rk(MW)&8=k=\sum s_i+\#\mbox{sections}&MW(\mathcal{E})\\
\hline
U\oplus E_8\oplus E_7&8+7& 4+3+1&\{1\}\\
\hline
U\oplus D_{14}\oplus A_1&14+1& 6+0+2&\Z/2\Z \\
\hline
\end{array}$\\

 $ $\\
\hline
\mbox{$g=2$,\ $k=7$, $r=16$,\ $a=2$,\ ($\delta$=1)}\\
$\begin{array}{|c|c|c|c|}
\hline
\mbox{trivial lattice}&14=\sum c_i+\rk(MW)&7=k=\sum s_i+\#\mbox{sections}&MW(\mathcal{E})\\
\hline
U\oplus E_7\oplus E_7&7+7& 3+3+1&\{1\}\\
\hline
U\oplus E_8\oplus D_6&8+6& 4+2+1&\{1\}\\
\hline
U\oplus D_{14}&14&6+1&\{1\}\\
\hline
U\oplus D_{12}\oplus A_1^2&12+1+1&5+0+0+2&\Z/2\Z\\
\hline
\end{array}$\\
 $ $\\
\hline
\mbox{$g=2$,\ $k=6$,\ $r=15$,\ $a=3$,\ ($\delta$=1)}\\
$\begin{array}{|c|c|c|c|}
\hline
\mbox{trivial lattice}&13=\sum c_i+\rk(MW)&6=k=\sum s_i+\#\mbox{sections}&MW(\mathcal{E})\\
\hline
U\oplus E_7\oplus D_6&7+6& 3+2+1&\{1\}\\
\hline
U\oplus E_8\oplus D_4\oplus A_1&8+4+1& 4+1+0+1&\{1\}\\
\hline
U\oplus D_{12}\oplus A_1&12+1&5+0+1&\{1\}\\
\hline
U\oplus D_{10}\oplus A_1^3&10+1+1+1&4+0+0+2&\Z/2\Z\\
\hline
\end{array}$\\
 $ $\\
\hline
\mbox{$g=2,\ k=5, r=14,\ a=4$, $\delta$=0 or 1}\\
$\begin{array}{|c|c|c|c|c|}
\hline
\mbox{trivial lattice}&12=\sum c_i+\rk(MW)&5=k=\sum s_i+\#\mbox{sections}&MW(\mathcal{E})&\delta\\
\hline
U\oplus D_6\oplus D_6&6+6& 2+2+1&\{1\}&1\\
\hline
U\oplus E_7\oplus D_4\oplus A_1&7+4+1& 3+1+0+1&\{1\}&1\\
\hline
U\oplus D_{10}\oplus A_1^2&10+1+1& 4+0+0+1&\{1\}&1\\
\hline
U\oplus E_8\oplus A_1^4&8+1+1+1+1&4+0+0+0+0+1&\{1\}&1\\
\hline
U\oplus D_8\oplus A_1^4&8+1+1+1+1&3+0+0+0+0+2&\Z/2\Z&1\\
\hline
U\oplus D_8\oplus D_4&8+4& 3+2+1&\{1\}&0\\
\hline
U\oplus E_7\oplus A_1^5&7+1+1+1+1+1&3+0+0+0+0+0+2&\Z/2\Z&0\\
\hline
\end{array}$\\
$ $\\
\hline
\mbox{$g=2, k=4, r=13,\ a=5$ ($\delta=1$)}\\
$\begin{array}{|c|c|c|c|}
\hline
\mbox{trivial lattice}&11=\sum c_i+\rk(MW)&4=k=\sum s_i+\#\mbox{sections}&MW(\mathcal{E})\\
\hline
U\oplus D_6\oplus D_4\oplus A_1&6+4+1&2+1+0+1&\{1\}\\
\hline
U\oplus D_8\oplus A_1^3&8+1+1+1&3+0+0+1&\{1\}\\
\hline
U\oplus E_7\oplus A_1^4&7+1+1+1+1&3+0+0+0+0+1&\{1\}\\
\hline
U\oplus D_6\oplus A_1^5&6+1+1+1+1+1&2+0+0+0+0+0+1&\Z/2\Z\\
\hline
\end{array}$\\
$ $\\
\hline
\mbox{$g=2,\ k=3, r=12,\ a=6$, ($\delta$=1)}\\
$\begin{array}{|c|c|c|c|}
\hline
\mbox{trivial lattice}&10=\sum c_i+\rk(MW)&3=k=\sum s_i+\#\mbox{sections}&MW(\mathcal{E})\\
\hline
U\oplus D_4\oplus D_4\oplus A_1^2&4+4+1+1&1+1+0+0+1&\{1\}\\
\hline
U\oplus D_6\oplus A_1^4&4+4+1+1+1+1&2+0+0+0+0+1&\{1\}\\
\hline
U\oplus D_4\oplus A_1^6&4+1+1+1+1+1+1&1+0+0+0+0+0+0+2&\Z/2\Z\\
\hline
\end{array}$\\
 $ $\\
\hline
\mbox{$g=2,\ k=2, r=11,\ a=7,$ ($\delta=1$)}\\
$\begin{array}{|c|c|c|c|}
\hline
\mbox{trivial lattice}&9=\sum c_i+\rk(MW)&2=k=\sum s_i+\#\mbox{sections}&MW(\mathcal{E})\\
\hline
U\oplus D_4\oplus A_1^5&4+1+1+1+1+1&1+0+0+0+0+0+1&\{1\}\\
\hline
U\oplus  A_1^9&1+1+1+1+1+1+1+1+1&0+0+0+0+0+0+0+0+0+2&\Z/2\Z\\
\hline
\end{array}$\\

$ $\\
\hline
\mbox{$g=2, k=1, r=10,\ a=8,$ ($\delta=1$)}\\
$\begin{array}{|c|c|c|c|}
\hline
\mbox{trivial lattice}&8=\sum c_i+\rk(MW)&2=k=\sum s_i+\#\mbox{sections}&MW(\mathcal{E})\\
\hline
U\oplus  A_1^8&1+1+1+1+1+1+1+1&0+0+0+0+0+0+0+0+1&\{1\}\\
\hline
\end{array}$
\end{longtable}

\begin{longtable}{l}\caption{}\label{table of admissible fibers g=3}\\ \hline $g=3, k=6, r=14, a=2$ ($\delta$=0)\\
$
\begin{array}{|c|c|c|c|c|}
\hline
\mbox{trivial lattice}&12=\sum c_i+\rk(MW)&6=k=\sum s_i+\#\mbox{sections}&MW(\mathcal{E})\\
\hline
U\oplus D_{12}&12&5+1&\{1\}\\
\hline
U\oplus E_8\oplus D_{4}&8+4& 4+1+1&\{1\}\\
\hline
\end{array}$\\

 $ $\\
\hline
\mbox{$g=3, k=5, r=13,\ a=3$, ($\delta$=1)}\\
$\begin{array}{|c|c|c|c|}
\hline
\mbox{trivial lattice}&11=\sum c_i+\rk(MW)&5=k=\sum s_i+\#\mbox{sections}&MW(\mathcal{E})\\
\hline
U\oplus E_7\oplus D_4&7+4& 3+1+1&\{1\}\\
\hline
U\oplus D_{10}\oplus A_1&10+1& 4+0+1&\{1\}\\
\hline
U\oplus E_8\oplus A_1^3&8+1+1+1&4+0+0+0+1&\{1\}\\
\hline
\end{array}$\\
 $ $\\
\hline
\mbox{$g=3, k=4, r=12,\ a=4$, ($\delta$=1)}\\
$\begin{array}{|c|c|c|c|}
\hline
\mbox{trivial lattice}&10=\sum c_i+\rk(MW)&4=k=\sum s_i+\#\mbox{sections}&MW(\mathcal{E})\\
\hline
U\oplus D_6\oplus D_4 &6+4&2+1+0+1&\{1\}\\
\hline
U\oplus D_8\oplus A_1^2&8+1+1&3+0+0+1&\{1\}\\
\hline
U\oplus E_7\oplus A_1^3&7+1+1+1&3+0+0+0+1&\{1\}\\
\hline
\end{array}$\\
 $ $\\
\hline
\mbox{$g=3, k=3, r=11,\ a=5,$ ($\delta$=1)}\\
$\begin{array}{|c|c|c|c|}
\hline
\mbox{trivial lattice}&9=\sum c_i+\rk(MW)&3=k=\sum s_i+\#\mbox{sections}&MW(\mathcal{E})\\
\hline
U\oplus D_4\oplus D_4\oplus A_1&4+4+1&1+1+0+1&\{1\}\\
\hline
U\oplus D_6\oplus A_1^3&4+4+1+1+1&2+0+0+0+1&\{1\}\\
\hline
\end{array}$\\
 $ $\\
\hline
\mbox{g=3, k=2, r=10,\ a=6, $\delta=$ 0 or 1}\\
$\begin{array}{|c|c|c|c|c|}
\hline
\mbox{trivial lattice}&8=\sum c_i+\rk(MW)&2=k=\sum s_i+\#\mbox{sections}&MW(\mathcal{E})&\delta\\
\hline
U\oplus D_4\oplus A_1^4&4+1+1+1+1&1+0+0+0+1&\{1\}&1\\
\hline
U\oplus  A_1^8&1+1+1+1+1+1+1&0+0+0+0+0+0+0+0+2&\Z/2\Z&0\\
\hline
\end{array}$\\

$ $\\
\hline
\mbox{$g=3, k=1, r=9,\ a=7,$ ($\delta=1$)}\\
$\begin{array}{|c|c|c|c|}
\hline
\mbox{trivial lattice}&8=\sum c_i+\rk(MW)&2=k=\sum s_i+\#\mbox{sections}&MW(\mathcal{E})\\
\hline
U\oplus  A_1^7&1+1+1+1+1+1+1&0+0+0+0+0+0+0+1&\{1\}\\
\hline
\end{array}$

\end{longtable}
\begin{longtable}{l}\caption{}\label{table of admissible fibers g=4}\\ \hline g=4,\ k=5, r=12, a=2, ($\delta=1$)\\
$
\begin{array}{|c|c|c|c|c|}
\hline
\mbox{trivial lattice}&10=\sum c_i+\rk(MW)&5=k=\sum s_i+\#\mbox{sections}&MW(\mathcal{E})\\
\hline
U\oplus D_{10}&10&4+1&\{1\}\\
\hline
U\oplus E_8\oplus A_1^2&8+1+1& 4+0+0+1&\{1\}\\
\hline
\end{array}$\\
$ $\\
\hline
\mbox{$g=4, k=4, r=11,\ a=3,$ ($\delta=1$)}\\
$\begin{array}{|c|c|c|c|}
\hline
\mbox{trivial lattice}&9=\sum c_i+\rk(MW)&4=k=\sum s_i+\#\mbox{sections}&MW(\mathcal{E})\\
\hline
U\oplus E_7\oplus A_1^2&7+1+1& 3+0+0+1&\{1\}\\
\hline
U\oplus D_{8}\oplus A_1&8+1& 3+0+1&\{1\}\\
\hline
\end{array}$\\
 $ $\\
\hline
\mbox{$g=4, k=3, r=10,\ a=4,$ $\delta=$ 0 or 1}\\
$\begin{array}{|c|c|c|c|c|}
\hline
\mbox{trivial lattice}&8=\sum c_i+\rk(MW)&3=k=\sum s_i+\#\mbox{sections}&MW(\mathcal{E})&\delta\\
\hline
U\oplus D_6\oplus A_1^2&6+1+1&2+0+0+1&\{1\}&1\\
\hline
U\oplus D_4\oplus D_4&4+4&1+1+0+1&\{1\}&0\\
\hline
\end{array}$\\
 $ $\\
\hline
\mbox{$g=4, k=2, r=9,\ a=5$, ($\delta=1$)}\\
$\begin{array}{|c|c|c|c|}
\hline
\mbox{trivial lattice}&7=\sum c_i+\rk(MW)&2=k=\sum s_i+\#\mbox{sections}&MW(\mathcal{E})\\
\hline
U\oplus D_4\oplus A_1^3&4+1+1+1&1+0+0+0+1&\{1\}\\
\hline
\end{array}$\\
$ $\\
\hline
\mbox{$g=4, k=1, r=8,\ a=6,$ ($\delta=1$)}\\
$\begin{array}{|c|c|c|c|}
\hline
\mbox{trivial lattice}&6=\sum c_i+\rk(MW)&1=k=\sum s_i+\#\mbox{sections}&MW(\mathcal{E})\\
\hline
U\oplus A_1^6&1+1+1+1+1+1&0+0+0+0+0+0+1&\{1\}\\
\hline
\end{array}$

\end{longtable}
\begin{longtable}{l}\caption{}\label{table of admissible fibers g=5}\\ \hline $g=5, k=5, r=11, a=1$, ($\delta=1$) \\
$
\begin{array}{|c|c|c|c|c|}
\hline
\mbox{trivial lattice}&9=\sum c_i+\rk(MW)&5=k=\sum s_i+\#\mbox{sections}&MW(\mathcal{E})\\
\hline
U\oplus E_8\oplus A_1&8+1& 4+0+1&\{1\}\\
\hline
\end{array}$\\
$ $\\
\hline
\mbox{$g=5, k=4, r=10,\ a=2$, $\delta$=0 or 1}\\
$\begin{array}{|c|c|c|c|c|}
\hline
\mbox{trivial lattice}&8=\sum c_i+\rk(MW)&4=k=\sum s_i+\#\mbox{sections}&MW(\mathcal{E})&\delta\\
\hline
U\oplus E_7\oplus A_1&7+1& 3+0+1&\{1\}&1\\
\hline
U\oplus D_{8}&8& 3+1&\{1\}&0\\
\hline
\end{array}$\\
 $ $\\
\hline
\mbox{$g=5, k=3, r=9,\ a=3$, ($\delta=1$)}\\
$\begin{array}{|c|c|c|c|}
\hline
\mbox{trivial lattice}&7=\sum c_i+\rk(MW)&3=k=\sum s_i+\#\mbox{sections}&MW(\mathcal{E})\\
\hline
U\oplus D_6\oplus A_1&6+1&2+0+1&\{1\}\\
\hline
\end{array}$\\
$ $\\
\hline
\mbox{$g=5, k=2, r=8,\ a=4$, ($\delta=1$)}\\
$\begin{array}{|c|c|c|c|}
\hline
\mbox{trivial lattice}&6=\sum c_i+\rk(MW)&2=k=\sum s_i+\#\mbox{sections}&MW(\mathcal{E})\\
\hline
U\oplus D_4\oplus A_1^2&4+1+1&1+0+0+1&\{1\}\\
\hline
\end{array}$\\
 $ $\\
\hline
\mbox{$g=5, k=1, r=7,\ a=5$}\\
$\begin{array}{|c|c|c|c|}
\hline
\mbox{trivial lattice}&5=\sum c_i+\rk(MW)&1=k=\sum s_i+\#\mbox{sections}&MW(\mathcal{E})\\
\hline
U\oplus A_1^5&1+1+1+1+1&0+0+0+0+0+1&\{1\}\\
\hline
\end{array}$
\end{longtable}
\begin{longtable}{l}\caption{}\label{table of admissible fibers g=6}\\ \hline $g=6, k=5, r=10, a=0$, ($\delta=0$)\\
$
\begin{array}{|c|c|c|c|c|}
\hline
\mbox{trivial lattice}&8=\sum c_i+\rk(MW)&5=k=\sum s_i+\#\mbox{sections}&MW(\mathcal{E})\\
\hline
U\oplus E_8&8& 4+1&\{1\}\\
\hline
\end{array}$\\
$ $\\
\hline
\mbox{$g=6, k=4, r=9,\ a=1$, ($\delta=1$)}\\
$\begin{array}{|c|c|c|c|}
\hline
\mbox{trivial lattice}&7=\sum c_i+\rk(MW)&4=k=\sum s_i+\#\mbox{sections}&MW(\mathcal{E})\\
\hline
U\oplus E_7&7& 3+1&\{1\}\\
\hline
\end{array}$\\
$ $\\
\hline
\mbox{$g=6, k=3, r=8,\ a=2$, ($\delta=1$)}\\
$\begin{array}{|c|c|c|c|}
\hline
\mbox{trivial lattice}&6=\sum c_i+\rk(MW)&3=k=\sum s_i+\#\mbox{sections}&MW(\mathcal{E})\\
\hline
U\oplus D_6&6&2+1&\{1\}\\
\hline
\end{array}$\\
 $ $\\
\hline
\mbox{$g=6, k=2, r=7,\ a=3$, ($\delta=1$)}\\
$\begin{array}{|c|c|c|c|}
\hline
\mbox{trivial lattice}&5=\sum c_i+\rk(MW)&2=k=\sum s_i+\#\mbox{sections}&MW(\mathcal{E})\\
\hline
U\oplus D_4\oplus A_1&4+1&1+0+1&\{1\}\\
\hline
\end{array}$\\
$ $\\
\hline
\mbox{$g=6, k=1, r=6,\ a=4$, ($\delta=1$)}\\
$\begin{array}{|c|c|c|c|}
\hline
\mbox{trivial lattice}&4=\sum c_i+\rk(MW)&1=k=\sum s_i+\#\mbox{sections}&MW(\mathcal{E})\\
\hline
U\oplus A_1^4&1+1+1+1&0+0+0+0+1&\{1\}\\
\hline
\end{array}$
\end{longtable}

\begin{longtable}{l}\caption{}\label{table of admissible fibers g=7}\\
\mbox{$g=7, k=2, r=6,\ a=2$, ($\delta=1$)}\\
$\begin{array}{|c|c|c|c|}
\hline
\mbox{trivial lattice}&4=\sum c_i+\rk(MW)&2=k=\sum s_i+\#\mbox{sections}&MW(\mathcal{E})\\
\hline
U\oplus D_4&4&1+1&\{1\}\\
\hline
\end{array}$\\
 $ $\\
\hline
\mbox{$g=7, k=1, r=5,\ a=3$, ($\delta=1$)}\\
$\begin{array}{|c|c|c|c|}
\hline
\mbox{trivial lattice}&4=\sum c_i+\rk(MW)&1=k= s_i+\#\mbox{sections}&MW(\mathcal{E})\\
\hline
U\oplus A_1^3&1+1+1&0+0+0+1&\{1\}\\
\hline
\end{array}$
\end{longtable}

\begin{longtable}{l}\caption{}\label{table of admissible fibers g=8+b}\\
\mbox{$g=8+b, k=1, r=4-b,\ a=2-b$, with $0\leq b\leq 2$ ($\delta$=1 if $b\leq 1$, $\delta=0$ if $b=2$)}\\
$\begin{array}{|c|c|c|c|}
\hline
\mbox{trivial lattice}&2-b=\sum c_i+\rk(MW)&2=k=\sum s_i+\#\mbox{sections}&MW(\mathcal{E})\\
\hline
U\oplus A_1^{2-b}&2-b&0+1&\{1\}\\
\hline
\end{array}$\\
\end{longtable}

{\it {\bf Acknowledgments}.
The authors would like to acknowledge the organization of Women in Numbers Europe I and II, where they met and started their collaboration. The authors would like to thank F. Balestrieri, J. Desjardins, C. Maistret, I. Vogt while the second author also thanks J. Top and M. Schuett for many interesting discussions.}

\end{document}